\catcode`\@=11

\font\teneufm=eufm10
\font\seveneufm=eufm7
\font\fiveeufm=eufm5
\newfam\eufmfam
\textfont\eufmfam=\teneufm  \scriptfont\eufmfam=\seveneufm
  \scriptscriptfont\eufmfam=\fiveeufm
\def\eufm@{\hexnumber@\eufmfam}

\font\tenmsa=msam10
\font\sevenmsa=msam7
\font\fivemsa=msam5
\font\tenmsb=msbm10
\font\sevenmsb=msbm7
\font\fivemsb=msbm5
\newfam\msafam
\newfam\msbfam
\textfont\msafam=\tenmsa  \scriptfont\msafam=\sevenmsa
  \scriptscriptfont\msafam=\fivemsa
\textfont\msbfam=\tenmsb  \scriptfont\msbfam=\sevenmsb
  \scriptscriptfont\msbfam=\fivemsb

\def\hexnumber@#1{\ifnum#1<10 \number#1\else
 \ifnum#1=10 A\else\ifnum#1=11 B\else\ifnum#1=12 C\else
 \ifnum#1=13 D\else\ifnum#1=14 E\else\ifnum#1=15 F\fi\fi\fi\fi\fi\fi\fi}

\def\msa@{\hexnumber@\msafam}
\def\msb@{\hexnumber@\msbfam}

\mathchardef\gx="2\eufm@78
\mathchardef\gg="2\eufm@67
\mathchardef\gm="2\eufm@6D
\mathchardef\gn="2\eufm@6E
\mathchardef\gd="2\eufm@64

\mathchardef\boxdot="2\msa@00
\mathchardef\boxplus="2\msa@01
\mathchardef\boxtimes="2\msa@02
\mathchardef\square="0\msa@03
\mathchardef\blacksquare="0\msa@04
\mathchardef\centerdot="2\msa@05
\mathchardef\lozenge="0\msa@06
\mathchardef\blacklozenge="0\msa@07
\mathchardef\circlearrowright="3\msa@08
\mathchardef\circlearrowleft="3\msa@09
\mathchardef\rightleftharpoons="3\msa@0A
\mathchardef\leftrightharpoons="3\msa@0B
\mathchardef\boxminus="2\msa@0C
\mathchardef\Vdash="3\msa@0D
\mathchardef\Vvdash="3\msa@0E
\mathchardef\vDash="3\msa@0F
\mathchardef\twoheadrightarrow="3\msa@10
\mathchardef\twoheadleftarrow="3\msa@11
\mathchardef\leftleftarrows="3\msa@12
\mathchardef\rightrightarrows="3\msa@13
\mathchardef\upuparrows="3\msa@14
\mathchardef\downdownarrows="3\msa@15
\mathchardef\upharpoonright="3\msa@16

\mathchardef\downharpoonright="3\msa@17
\mathchardef\upharpoonleft="3\msa@18
\mathchardef\downharpoonleft="3\msa@19
\mathchardef\rightarrowtail="3\msa@1A
\mathchardef\leftarrowtail="3\msa@1B
\mathchardef\leftrightarrows="3\msa@1C
\mathchardef\rightleftarrows="3\msa@1D
\mathchardef\Lsh="3\msa@1E
\mathchardef\Rsh="3\msa@1F
\mathchardef\rightsquigarrow="3\msa@20
\mathchardef\leftrightsquigarrow="3\msa@21
\mathchardef\looparrowleft="3\msa@22
\mathchardef\looparrowright="3\msa@23
\mathchardef\circeq="3\msa@24
\mathchardef\succsim="3\msa@25
\mathchardef\gtrsim="3\msa@26
\mathchardef\gtrapprox="3\msa@27
\mathchardef\multimap="3\msa@28
\mathchardef\therefore="3\msa@29
\mathchardef\because="3\msa@2A
\mathchardef\doteqdot="3\msa@2B

\mathchardef\triangleq="3\msa@2C
\mathchardef\precsim="3\msa@2D
\mathchardef\lesssim="3\msa@2E
\mathchardef\lessapprox="3\msa@2F
\mathchardef\eqslantless="3\msa@30
\mathchardef\eqslantgtr="3\msa@31
\mathchardef\curlyeqprec="3\msa@32
\mathchardef\curlyeqsucc="3\msa@33
\mathchardef\preccurlyeq="3\msa@34
\mathchardef\leqq="3\msa@35
\mathchardef\leqslant="3\msa@36
\mathchardef\lessgtr="3\msa@37
\mathchardef\backprime="0\msa@38
\mathchardef\risingdotseq="3\msa@3A
\mathchardef\fallingdotseq="3\msa@3B
\mathchardef\succcurlyeq="3\msa@3C
\mathchardef\geqq="3\msa@3D
\mathchardef\geqslant="3\msa@3E
\mathchardef\gtrless="3\msa@3F
\mathchardef\sqsubset="3\msa@40
\mathchardef\sqsupset="3\msa@41
\mathchardef\trianglerighteq="3\msa@44
\mathchardef\trianglelefteq="3\msa@45
\mathchardef\bigstar="0\msa@46
\mathchardef\between="3\msa@47
\mathchardef\blacktriangledown="0\msa@48
\mathchardef\blacktriangleright="3\msa@49
\mathchardef\blacktriangleleft="3\msa@4A
\mathchardef\blacktriangle="0\msa@4E
\mathchardef\triangledown="0\msa@4F
\mathchardef\eqcirc="3\msa@50
\mathchardef\lesseqgtr="3\msa@51
\mathchardef\gtreqless="3\msa@52
\mathchardef\lesseqqgtr="3\msa@53
\mathchardef\gtreqqless="3\msa@54
\mathchardef\Rrightarrow="3\msa@56
\mathchardef\Lleftarrow="3\msa@57
\mathchardef\veebar="2\msa@59
\mathchardef\barwedge="2\msa@5A
\mathchardef\doublebarwedge="2\msa@5B
\mathchardef\angle="0\msa@5C
\mathchardef\measuredangle="0\msa@5D
\mathchardef\sphericalangle="0\msa@5E
\mathchardef\varpropto="3\msa@5F
\mathchardef\smallsmile="3\msa@60
\mathchardef\smallfrown="3\msa@61
\mathchardef\Subset="3\msa@62
\mathchardef\Supset="3\msa@63
\mathchardef\Cup="2\msa@64

\mathchardef\Cap="2\msa@65

\mathchardef\curlywedge="2\msa@66
\mathchardef\curlyvee="2\msa@67
\mathchardef\leftthreetimes="2\msa@68
\mathchardef\rightthreetimes="2\msa@69
\mathchardef\subseteqq="3\msa@6A
\mathchardef\supseteqq="3\msa@6B
\mathchardef\bumpeq="3\msa@6C
\mathchardef\Bumpeq="3\msa@6D
\mathchardef\lll="3\msa@6E

\mathchardef\ggg="3\msa@6F

\mathchardef\circledS="0\msa@73
\mathchardef\pitchfork="3\msa@74
\mathchardef\dotplus="2\msa@75
\mathchardef\backsim="3\msa@76
\mathchardef\backsimeq="3\msa@77
\mathchardef\complement="0\msa@7B
\mathchardef\intercal="2\msa@7C
\mathchardef\circledcirc="2\msa@7D
\mathchardef\circledast="2\msa@7E
\mathchardef\circleddash="2\msa@7F
\def\ulcorner{\delimiter"4\msa@70\msa@70 }
\def\urcorner{\delimiter"5\msa@71\msa@71 }
\def\llcorner{\delimiter"4\msa@78\msa@78 }
\def\lrcorner{\delimiter"5\msa@79\msa@79 }
\def\yen{\mathhexbox\msa@55 }
\def\checkmark{\mathhexbox\msa@58 }
\def\circledR{\mathhexbox\msa@72 }
\def\maltese{\mathhexbox\msa@7A }
\mathchardef\lvertneqq="3\msb@00
\mathchardef\gvertneqq="3\msb@01
\mathchardef\nleq="3\msb@02
\mathchardef\ngeq="3\msb@03
\mathchardef\nless="3\msb@04
\mathchardef\ngtr="3\msb@05
\mathchardef\nprec="3\msb@06
\mathchardef\nsucc="3\msb@07
\mathchardef\lneqq="3\msb@08
\mathchardef\gneqq="3\msb@09
\mathchardef\nleqslant="3\msb@0A
\mathchardef\ngeqslant="3\msb@0B
\mathchardef\lneq="3\msb@0C
\mathchardef\gneq="3\msb@0D
\mathchardef\npreceq="3\msb@0E
\mathchardef\nsucceq="3\msb@0F
\mathchardef\precnsim="3\msb@10
\mathchardef\succnsim="3\msb@11
\mathchardef\lnsim="3\msb@12
\mathchardef\gnsim="3\msb@13
\mathchardef\nleqq="3\msb@14
\mathchardef\ngeqq="3\msb@15
\mathchardef\precneqq="3\msb@16
\mathchardef\succneqq="3\msb@17
\mathchardef\precnapprox="3\msb@18
\mathchardef\succnapprox="3\msb@19
\mathchardef\lnapprox="3\msb@1A
\mathchardef\gnapprox="3\msb@1B
\mathchardef\nsim="3\msb@1C
\mathchardef\napprox="3\msb@1D
\mathchardef\nsubseteqq="3\msb@22
\mathchardef\nsupseteqq="3\msb@23
\mathchardef\subsetneqq="3\msb@24
\mathchardef\supsetneqq="3\msb@25
\mathchardef\subsetneq="3\msb@28
\mathchardef\supsetneq="3\msb@29
\mathchardef\nsubseteq="3\msb@2A
\mathchardef\nsupseteq="3\msb@2B
\mathchardef\nparallel="3\msb@2C
\mathchardef\nmid="3\msb@2D
\mathchardef\nshortmid="3\msb@2E
\mathchardef\nshortparallel="3\msb@2F
\mathchardef\nvdash="3\msb@30
\mathchardef\nVdash="3\msb@31
\mathchardef\nvDash="3\msb@32
\mathchardef\nVDash="3\msb@33
\mathchardef\ntrianglerighteq="3\msb@34
\mathchardef\ntrianglelefteq="3\msb@35
\mathchardef\ntriangleleft="3\msb@36
\mathchardef\ntriangleright="3\msb@37
\mathchardef\nleftarrow="3\msb@38
\mathchardef\nrightarrow="3\msb@39
\mathchardef\nLeftarrow="3\msb@3A
\mathchardef\nRightarrow="3\msb@3B
\mathchardef\nLeftrightarrow="3\msb@3C
\mathchardef\nleftrightarrow="3\msb@3D
\mathchardef\divideontimes="2\msb@3E
\mathchardef\varnothing="0\msb@3F
\mathchardef\nexists="0\msb@40
\mathchardef\mho="0\msb@66
\mathchardef\thorn="0\msb@67
\mathchardef\beth="0\msb@69
\mathchardef\gimel="0\msb@6A
\mathchardef\daleth="0\msb@6B
\mathchardef\lessdot="3\msb@6C
\mathchardef\gtrdot="3\msb@6D
\mathchardef\ltimes="2\msb@6E
\mathchardef\rtimes="2\msb@6F
\mathchardef\shortmid="3\msb@70
\mathchardef\shortparallel="3\msb@71
\mathchardef\smallsetminus="2\msb@72
\mathchardef\thicksim="3\msb@73
\mathchardef\thickapprox="3\msb@74
\mathchardef\approxeq="3\msb@75
\mathchardef\succapprox="3\msb@76
\mathchardef\precapprox="3\msb@77
\mathchardef\curvearrowleft="3\msb@78
\mathchardef\curvearrowright="3\msb@79
\mathchardef\digamma="0\msb@7A
\mathchardef\varkappa="0\msb@7B
\mathchardef\hslash="0\msb@7D
\mathchardef\hbar="0\msb@7E
\mathchardef\backepsilon="3\msb@7F
\def\Bbb{\ifmmode\let\next\Bbb@\else
 \def\next{\errmessage{Use \string\Bbb\space only in math mode}}\fi\next}
\def\Bbb@#1{{\Bbb@@{#1}}}
\def\Bbb@@#1{\fam\msbfam#1}

\mathchardef\cg="2\eufm@67
\mathchardef\cm="2\eufm@6D

\catcode`\@=\active

\def\CL{\hbox{{$\cal L$}}}

  \def\CE{\hbox{{$\cal E$}}}

\def\R{{\Bbb R}}
\def\C{{\Bbb C}}
\def\Z{{\Bbb Z}}

\def\h{{{1\over2}}}
\def\eps{{\epsilon}}

\def\trace{{\rm Tr\, }}

\def\<{\langle}
\def\>{\rangle}

\def\lcross{{>\!\!\!\triangleleft}}

\def\cobicross{{\triangleright\!\!\!\blacktriangleleft}}
\def\bicross{{\blacktriangleright\!\!\!\triangleleft}}
\def\dcross{{\bowtie}}
\def\codcross{{\blacktriangleright\!\!\blacktriangleleft}}

\def\tens{\mathop{\otimes}}
\def\la{{\triangleright}}\def\ra{{\triangleleft}}
\def\isom{{\cong}}

\def\Ad{{\rm Ad}}
\def\ad{{\rm ad}}

\def\pisl{{\pi\!\!\! /}}
\def\id{{\rm id}}

\def\extd{{\rm d}}

\def\text#1{{\rm #1}}
\def\note#1{}

\def\nquad{{\!\!\!\!\!\!}}

\def\equad{\nquad}
\def\eqn#1#2{\begin{equation}#2\label{#1}\end{equation}}
\def\align#1{\begin{eqnarray*}#1\end{eqnarray*}}
\documentstyle[11pt]{article}
\begin{document}
\baselineskip 22pt

{\ }\qquad  \hskip 4in Damtp/1999-72
\vspace{.2in}

\begin{center} {\LARGE POISSON-LIE T-DUALITY FOR QUASITRIANGULAR LIE
BIALGEBRAS}
\\ \baselineskip 13pt{\ }\\
{\ }\\ Shahn Majid\footnote{Royal Society University Research
Fellow and Fellow of Pembroke College, Cambridge}\\ {\ }\\
Department of Applied Mathematics \& Theoretical Physics\\
University of Cambridge, Cambridge CB3 9EW, UK\\ +\\ E.J.
Beggs\\{\ }\\Department of Mathematics,\\University of Wales,
\\  Swansea SA2 8PP, UK
\end{center}
\begin{center}
May, 1999
\end{center}
%\vspace{10pt}

\begin{quote}\baselineskip 13pt
\noindent{\bf Abstract}  We introduce a new 2-parameter family of sigma
models exhibiting Poisson-Lie T-duality on a quasitriangular
Poisson-Lie group $G$. The models contain previously known models
as well as a new 1-parameter line of models having the novel feature that the
Lagrangian takes the simple form ${\cal L}=E(u^{-1}u_+,u^{-1}u_-)$
where the generalised metric $E$ is constant ({\em not} dependent on
the field $u$ as in previous models). We characterise these models in 
terms of a global
conserved $G$-invariance. The models on $G=SU_2$ and its dual
$G^\star$ are
computed explicitly. The general theory of Poisson-Lie
T-duality is also extended; we develop the Hamiltonian
formulation and the reduction for constant loops to
integrable motion on the group manifold. Finally, we
generalise T-duality in the Hamiltonian formulation to group
factorisations $D=G\dcross M$ where the
subgroups need not be dual or even have the same dimension and need
not be connected to the Drinfeld double or to Poisson structures.

\end{quote}
\baselineskip 22pt
\section{Introduction}

Poisson-Lie T-duality has been introduced in \cite{Kli:poi}\cite{KliSev:poi}
and other works
as a non-Abelian version of T-duality in string theory, based on duality of
Lie bialgebras.  A motivation (stated in
\cite{Kli:poi}) is quantum group or Hopf algebra duality;  this had been
introduced as a duality for physics
several years previously\cite{Ma:pla}\cite{Ma:the}\cite{Ma:phy}\cite{Ma:mat},
as an `observable-state'
duality for certain quantum systems based on group factorisations 
$D=G\dcross M$.
In one system a particle moves
in $G$ under the action of $M$ and its quantum algebra of observables
 is the bicrossproduct Hopf algebra
$U(\gm)\cobicross \C(G)$, in the dual system the roles of $G,M$ are
interchanged
but its quantum algebra of observables $\C(M)\bicross U(\gg)$
has the same physical content
but with the roles of observables/states and position/momentum
interchanged (here $\gg,\gm$ are the Lie algebras of $G,M$ respectively). 
Indeed, being mutually
dual Hopf algebras the two quantum systems are related to each other by quantum
Fourier transform
\eqn{fou}{ {\cal F}: U(\gm)\cobicross \C(G)\to \C(M)\bicross U(\gg),}
see \cite{MaOec:dif} where this was recently studied in detail for the simplest
example (the so-called Planck-scale Hopf algebra $\C[p]\cobicross\C[x]$ in
\cite{Ma:pla}.)  Under this observable-state duality it was shown in 
\cite{Ma:pla} that
one had inversion of coupling constants as well as connections with 
Planck-scale
physics. At about the same time, Abelian T-duality was introduced in 
\cite{Tse:dua} and
elsewhere as a momentum-winding mode symmetry in string theory with some 
similar
features. The observable-state duality (\ref{fou}) is not, however, limited
in any way to the Abelian case and
indeed there is a natural model for every
compact simple group $G$ with $M=G^\star$
the Yang-Baxter dual. Here a Lie bialgebra is an infinitesimal version of a 
Hopf
algebra and has a dual $\gg^\star$, and $G^\star$ is its associated 
Lie group. It is also the group
of dressing transformations\cite{Sem:dre} in the theory of classical inverse 
scattering
and the solvable group in the Isawasa decomposition $D=G_{\C}
=G\dcross G^\star$ of
the complexification of the compact Lie group $G$, see \cite{Ma:mat}. Moreover,
$D=G\dcross G^\star$
is the Lie group associated to the Drinfeld double $\gd(\gg)$ of $\gg$ as
a Lie bialgebra \cite{Dri:ham}. The Lie bialgebra structure of $\gg$ also
implies a natural Poisson bracket on $G$\cite{Dri:ham}. Further details are in
the Preliminaries; see also \cite{Ma:book} for an introduction to these 
topics. These
 quantum systems $U(\gg)\cobicross \C(G^\star)$ with observable-state duality
 were constructed in \cite{Ma:the}\cite{Ma:phy}\cite{Ma:mat} as
one of the two main sources of quantum groups canonically associated to
a simple Lie algebra (the other is the more well-known q-deformation of 
$U(\gg)$ to quantum groups $U_q(\gg)$).

The subsequent theory of Poisson-Lie T-duality\cite{KliSev:poi} indeed has many
of the same
features. One system
consists of a sigma model on the group $G$ with a Lagrangian of the form
\[ \CL=E_{u}(u^{-1}u_+,u^{-1}u_-),\quad u:\R^{1,1}\to G,\]
where $u$ is the field, $u_\pm$ are derivatives in light-cone
coordinates and $E_u$ a bilinear form on $\gg$ but depending on the
value of $u$ (a `generalised metric' since $E_u$ need not be symmetric). The
dual theory is a sigma-model on $G^\star$ with
\[ \hat{\cal L}=\hat E_t(t^{-1}t_+,t^{-1}t_-),\quad t:\R^{1,1}\to G^\star.\]
The physical content of the two theories is established to be the same due
to the existence of the larger group $D=G\dcross G^\star$ associated to
 the Drinfeld double $\gd(\gg)$.

In the present paper we extend Poisson-Lie T-duality in
several directions, motivated in part by the above connections with quantum
groups and observable-state duality. From a physical point of view the
main result is as follows:
the previously-known models exhibiting Poisson-Lie T-duality require a very 
special
form of the generalised metric $E_u$ depending on $u$ in a
rather complicated way (related to the Poisson bracket on $G$). This is in 
sharp
contrast to the usual principal
sigma model\cite{Nov}  where the metric is a constant, the Killing form $K$. 
As a
result, Poisson-Lie T-duality would appear to be somewhat
artificial and to apply to only certain highly non-linear models
where the `metric' in the target group is far from constant.
Even the explicit form of $E_u$
is known only in some simple cases such as $\gg=b_+$ the Borel-subalgebra
of $su_2$ \cite{Kli:poi}. The $\gg=su_2$ case was discussed recently
in \cite{LleVar:poi} but still without fully explicit formulae for the 
resulting
Lagrangians. Our main result is the introduction of a new 2-parameter class
of models within the existing general framework for Poisson-Lie T-duality
but which much nicer
properties. We also provide new computational tools using the theory of Lie
bialgebras to compute the models explicitly. We obtain, for example, the 
explicit
Lagrangians in the $SU_2$ case and its dual.

These new models require that
$\gg$ is a quasitriangular Lie
bialgebra, i.e. defined by an element
$r\in\gg\tens\gg$ obeying the so-called modified {\em classical
Yang-Baxter equations}\cite{Dri:ham}. This includes
all complex semisimple Lie algebras equipped, for example, with
their standard Drinfeld-Sklyanin quasitriangular structure as used
in the theory of classical inverse scattering. The quantisations
of the associated Poisson bracket on $G$ in these cases include coordinate
algebras of the quantum groups $U_q(\gg)$. This is
therefore an important class of models, and we will find quite
tractable formulae in this case. We use $r$ not only in the Lie
bialgebra structure (which is usual) but again in certain boundary
conditions for the graph coordinates in order to cancel their natural
$u$-dependence for the choice of certain parameters. This greater
generality allows for a two-parameter family of models associated
to this data. Moreover, in this extended parameter space there is
a novel line of `nice' models in which $E_u=E_e$ is a {\em constant}
not dependent at all on $u$. This line includes at $\infty$ the standard
principal sigma model where $E_e=K$ the Killing form, but at other points
has an antisymmetric part built from $r$ itself. In this way one
may approach the principal sigma model itself along a line of sigma models
exhibiting Poisson-Lie T-duality and of a simple form without additional
non-linearities due to a non-constant generalised metric. The dual models
are more complicated but at $\infty$, for example, one obtains an Abelian 
model as the
Poisson-Lie T-dual of the principal sigma model approached in this way
(the latter lies on the boundary of the space of models exhibiting T-duality).
These results are presented in Section~6.

Also in the paper we develop the Hamiltonian picture of Poisson-Lie T-duality
in rather more detail than we have found elsewhere; see also \cite{Sfe:can}. 
This is
done in Section~3 after the preliminary Section~2. Among the new
results is a more regular expression for the Hamiltonian that covers both
the model and the dual model simultaneously. Also new is a study
of the symmetries of the theory induced by the {\em left} action of $D$ on
itself. These are not usually considered because they are not conserved
but we show that they do respect the symplectic structure. Moreover, when
$E_u$ is constant we show that the action of $G\subset D$ {\em is} conserved
and we compute the conserved charges.

A second general development, in Section~4, is a study of the classical 
mechanical
system on $G$ (say) in the limit of point-like strings (i.e. $x$-independent
solutions). We show that this constraint commutes with the dynamics and
 we provide the resulting Lagrangian and Hamiltonian systems and the phase 
space.
The left action of $D$ descends to the classical mechanical
system and we  show that it has a moment map. The conserved charges are
computed in the case of constant $E_u$.
The dual model on $G^\star$ equivalent to these point-solutions are not point
solutions but extended solutions of a certain special form.
We also discuss the quantisation of this classical mechanical
system both conventionally and in a manner relevant to the conserved
charges.  Although
these  systems appear to be different from the systems
$U(\gg^\star)\cobicross\C[G]$
exhibiting observable-state duality at the Planck-scale\cite{Ma:pla}, we do
establish some points of comparison, such as a common phase space.

Section~5 contains some further algebraic preliminaries needed
for the explicit construction of $E_u$. We show that
\[ \Ad^*_u(E_u)=(E_e^{-1}+\Pi(u))^{-1}\]
where $\Pi$ is the $\gg\tens\gg$-valued function
defining the Poisson-structure on $G$. To our knowledge this derivation differs
from previous work in that we do not assume anything about $E_e^{-1}$,
in particular it need not be the Killing form usually
added\cite{KliSev:poi} to $\Pi$ as an ansatz. This greater generality allows us
in Section~6 to present our main result; the class of `nice'
Poisson-Lie T-dual models based on quasitriangular Lie bialgebras.

Finally, Section~7 introduces new `double-Neumann' boundary
conditions for the open string and proceeds for these (as well
as more trivially for closed strings) to extend the
Poisson-Lie T-duality in the Hamiltonian form
to general group factorisations $D=G\dcross M$,
where $D$ need no longer be the Lie group of the Drinfeld double $\gd(\gg)$
and indeed $\gm$ need not be $\gg^\star$ but could be some quite
different Lie algebra, possibly of different dimension. This is directly 
motivated
by the observable-state duality models which exist\cite{Ma:the}\cite{Ma:book}
 for any factorisation. It is also motivated by the Adler-Kostant-Symes theorem
in classical inverse scattering which works for a general factorisation
equipped with an inner product, see \cite{Ma:book}.
 The dynamics are determined, similarly to the
conventional bialgebra theory, by the splitting of the Lie algebra
of $D$ into orthogonal subspaces but these need no longer be of
the same dimension (although only in this case is there a
sigma-model interpretation). We also have an action of $D$ by left
multiplication on the phase space with the double-Neumann boundary
conditions which us useful even for standard Poisson-Lie T-duality
based on Lie bialgebras. In particular, it extends to an action of
the affine Kac-Moody Lie algebra $\tilde \gd$.

Several directions remain for further work. First of all,
 only some first steps are taken (in Section~4)
to relate T-duality to observable-state duality (\ref{fou})
in the quantum theory; our long term motivation here is to extend these
ideas from particles to loops and hence to formulate T-duality for the full
quantum systems  as a duality operation on a more general algebraic
structure (no doubt more general than
Hopf algebras but in the same spirit). This in turn would give insight into
the correct algebraic structure for the conjectured
`M-theory' about which little is known  beyond dualities visible
in the Lagrangians at various classical limits. Let us mention only
that Poisson-Lie T-duality is connected also
with mirror symmetry\cite{Par:mir} and indirectly with several other
relevant dualities in the theory of strings and branes.

Secondly, there are some interesting examples of
the generalisation of Poisson-Lie T-duality in Section~7 which exist
in principle and should be developed further.
Thus, the conformal group on $\R^n$ ($n>2$) has, locally,
a factorisation into the Poincar\'e group and an $\R^n$ of special conformal
translations. The global structure of the factorisation is singular in a
similar manner to the `black-hole event-horizon'-like features
of the Planck-scale Hopf algebra $\C[p]\cobicross\C[x]$ in \cite{Ma:pla}.
There is also the possibility in our more general setting of a
many-sided T-duality (i.e. not only two equivalent theories) associated
to more than one factorisation of the same
group.

Finally, the natural emergence of generalised metrics
which have both symmetric and antisymmetric parts is a natural
feature of noncommutative Riemannian geometry\cite{Ma:rie} (where symmetry 
is natural
only in the commutative limit). This is a further
direction that remains to be explored. Also to be considered is
the addition of WZNW terms to render our 2-parameter class of
sigma-models conformally invariant as well as the computation of 1-loop or 
higher
quantum effects c.f. \cite{Sfe:poi} \cite{AKT:poi}.

\subsection*{Preliminaries}

We recall, see e.g.\cite{Ma:book} that a Lie bialgebra is a Lie algebra
equipped with $\delta:\gg\to \gg\tens\gg$ where $\delta$ is
antisymmetric and obeys the coJacobi identity (so that $\gg^*$ is
a Lie algebra) and
\[ \delta[\xi,\eta]=\ad_\xi(\eta)-\ad_\eta(\xi)\]
for all $\xi,\eta\in \gg$, where $\ad$ extends as a derivation.

Next, associated to any Lie bialgebra $\gg$ there is a double Lie
algebra $\gd=\gg\dcross\gg^{*{\rm op}}$. This is a double
semidirect sum with cross relations
\[ [\phi,\xi]=\phi\la\xi-\phi\ra\xi\] where the actions are
mutually coadjoint ones
\[ \phi\la\xi=\<\xi_{[2]},\phi\>\xi_{[1]},\quad
\phi\ra\xi=\<\xi,\phi_{[2]}\>\phi_{[1]}\]
where the angle brackets are the dual pairing of $\gg^{*}$ with
$\gg$ and $\delta(\xi)=\xi_{[1]}\tens \xi_{[2]}$. Here $\gd$ is
quasitriangular and factorisable (see later) and as a result there
is an adjoint invariant inner product on $\gd$,
$$
(\xi\oplus\phi,\eta\oplus\psi)\ =\ \< \phi,\eta\>\ +\
\<\phi,\xi\>\ .
$$
Here
\eqn{gstar}{ \gg^\star=\gg^{*\rm op}}
and $\gg$ are maximal isotropic subspaces. We will need this
description from \cite{Ma:phy} which is somewhat more explicit than the usual
description in terms of `Manin triple' in Drinfeld's work
\cite{Dri}.

Given a double cross sum of Lie algebras $\gg\dcross\gm$, we may
at least locally exponentiate to a double cross product of Lie
groups $G\dcross M$. This is given explicitly in \cite{Ma:mat}. We view the
Lie algebra actions as cocycles, exponentiate to Lie group
cocycles, view these as flat connections and take the parallel
transport operation. The actions can be described by $b(u)\in
\gg\tens \gm^*$ given by $b(u)(\phi)=b_\phi(u)=(\phi\la u)u^{-1}$
and $a(s)\in \gg^*\tens m$ given by
$a(s)(\xi)=a_\xi(s)=s^{-1}(s\ra\xi)$.
 It can be shown that  $b\in Z^1_{{\rm Ad}\tens \ra^*}(G,\gg\tens\gm^*)$
is a cocycle, where the action $\ra^*$ is a left action of $G$ on
$\gm^*$ given by dualising the right action $\ra:\gm\times G\to
\gm$. Also $a\in Z^1_{\la^*\tens {\rm Ad}_R}(M,\gg^*\tens\gm)$,
where ${\rm Ad}_R$ is the right adjoint action of $M$ on $\gm$ and
$\la^*$ is the right action of $M$ on $\gg^*$ given by dualising
its action on $\gg$. These Lie-algebra-valued functions $a,b$
generate the vector fields for the action of $\gg$ on $M$ and
$\gm$ on $G$ respectively. Thus, $\phi\la u=b_\phi(u)u$ where $\xi
u=\tilde{\xi}$ denotes the right invariant vector field on $G$
generated by $\xi\in\gg$. Similarly, $s\ra\xi=sa_\xi(s)$. Once the
global actions of $G$ on $M$ and vice-versa are known, the
structure of $G\dcross M$ is such that
\eqn{dcross}{su=(s\la u)(u\ra s),\quad \forall u\in G,\ s\in
M.} This allows every element of the double cross product group
$G\dcross M$ to be uniquely factorised either as $GM$ or as $MG$,
and relates the two factorisations.

\section{ T-Duality based on Lie bialgebras}

We begin by giving a version of the standard T-duality based on the
Drinfeld double of a Lie bialgebra \cite{Kli:poi}\cite{KliSev:poi}. We will
phrase it slightly
differently in terms of double cross products with a view to later
generalisation. Thus, there is a double cross product group
$D=G\dcross M$ with Lie algebra $\gd=\gg+\gm$, and an
adjoint-invariant bilinear form on $\gd$ which is zero on
restriction to $\gg$ and $\gm$. The Lie algebra $\gd$ is the direct
sum of two perpendicular
 subspaces $\cal E_-$ and $\cal E_+$.
This means that $\gm=\gg^{*op}$, that the factorisation is a
coadjoint matched pair and that $\gd=D(\gg)$, the Drinfeld double
of $\gg$, which is the setting that Klim\v{c}\'{\i}k etc., assume.

On $\Bbb R^2$ we use light cone coordinates $x_+=t+x$ and $x_-=t-x$,
where $t$ and $x$ are the standard time-space coordinates.
Now let us suppose that there is a function $k:\Bbb
R^2\to G\dcross M$, with the properties that $k_+k^{-1}(x_+,x_-)\in
\cal E_-$
and $k_-k^{-1}(x_+,x_-)\in \cal E_+$ for all $(x_+,x_-)\in \Bbb R^2$.
Then
we see that,
if we factor $k=us$ for $u\in G$ and $s\in M$,
\[
u^{-1}u_\pm\ +\ s_\pm s^{-1}\in u^{-1}{\cal E}_\mp u\ .
\]
If the projection $\pi_{\gg}:\gd\to \gg$ (with kernel $\gm$)
 is 1-1 and onto when restricted to
$u^{-1}{\cal E}_-u$ and
 $u^{-1}{\cal E}_+u$,
we can find graph coordinates $E_u:\gg\to \gm$ and $T_u:\gg\to
\gm$ so that
\[
\big\{ \xi+E_u(\xi):\xi\in\gg\big\}\ =\ u^{-1} {\cal E}_+u
\quad{\rm
and}\quad
\big\{ \xi+T_u(\xi):\xi\in\gg\big\}\ =\ u^{-1} {\cal E}_-u\ .
\]
It follows that $s_-s^{-1}=E_u(u^{-1}u_-)$ and
$s_+s^{-1}=T_u(u^{-1}u_+)$. From the identity
\[
(s_+s^{-1})_--(s_-s^{-1})_+=[s_-s^{-1},s_+s^{-1}]
\]
we deduce that $u(x_+,x_-)$ satisfies the equation
\eqn{ueqm}{
\big(T_u(u^{-1}u_+)\big)_-\ -\ \big(E_u(u^{-1}u_-)\big)_+\ =\
\big[E_u(u^{-1}u_-),T_u(u^{-1}u_+)\big]\ .}
Klim\v{c}\'{\i}k shows that the Lagrangian density
\eqn{ulagr}{
{\cal L}\ =\ \< E_u(u^{-1}u_-),u^{-1}u_+\>}
gives rise to these equations of motion.

The dual theory is given by the factorisation $k=tv$, where $t\in M$ and
$v\in G$. If we let
$\hat E_t:\gm\to\gg$ and $\hat T_t:\gm\to\gg$ be the graph coordinates
of $t^{-1}{\cal E}_+t$ and $t^{-1}{\cal E}_-t$ respectively, then $t(x_+,x_-)$
obeys the dual equation
\eqn{teqm}{
\big(\hat T_t(t^{-1}t_+)\big)_-\ -\ \big(\hat E_t(t^{-1}t_-)\big)_+\
=\ \big[\hat E_t(t^{-1}t_-),\hat T_t(t^{-1}t_+)\big]\ .}
These are the equations of motion for a sigma model with Lagrangian
\eqn{tlagr}{
\hat{\cal L}\ =\ \< \hat E_t(t^{-1}t_-),t^{-1}t_+\>.}
These two models are different but equivalent descriptions of the
model defined by $k$. The $(u,s)$ and $(t,v)$ coordinates are
related by the actions of the double cross product group
structure:
\eqn{utreln}{ tv=(t\la v)(t\ra v)=us.}

\section{Hamiltonian formulation of T-duality}
  There are two models considered in the last section, the first order
equations of motion for
$k:\Bbb R^2\to G\dcross M$ and the second order equations of motion for
$u:\Bbb R^2\to G$. The equations of motion for $k:\Bbb R^2\to G\dcross M$
are the natural way to introduce duality into the system, and are very nearly
equivalent to the equations of motion for  $u:\Bbb R^2\to G$. There is not
a 1-1
correspondence between the systems, as multiplying $k$ on the right by a
constant element of $M$ gives rise to exactly the same $u$. We have a
Lagrangian and Hamiltonian
for the $u$ equations of motion, and can work out the corresponding Hamiltonian
mechanics. However the reader must remember that this will not give the
Hamiltonian mechanics for $k$, but rather for $k$ quotiented on the right by
constant elements of $M$.

As pointed out by Klim\v{c}\'{\i}k, we
can take the phase space of the system to be the set of smooth
functions $C^\infty(\Bbb R,D)$
(or more strictly $C^\infty(\Bbb R,D)/M$), where we regard $\Bbb R$ to be a
constant time line in $\Bbb R^{1+1}$, or
$C^\infty\big((0,\pi),D\big)/M$ for a finite space. We will compute
the symplectic structure more explicitly than we have found
elsewhere and then obtain a new and more symmetric formulation of
the Hamiltonian density that covers both the model and the dual
model simultaneously. We will need this in later sections when we
generalise to arbitrary factorisations, as well as for the point-like
limit.

\subsection{The symplectic form}

We begin by showing that this is the correct phase space, i.e. that
such a function encodes both $u$ and $\dot u$ on a constant time
line. Thus, take $k\in C^\infty(\Bbb R,D)$ or
$C^\infty\big((0,\pi),D\big)$. As $k(x)\in D$ we can factor it as
$k(x)=u(x)s(x)$, so $u(x)$ is specified on the constant time line.
But we also know that
\eqn{sx}{
s_xs^{-1}\ =\ T_u(u^{-1}u_+)-E_u(u^{-1}u_-)\ =\ \frac12\Big(
T_u(u^{-1}\dot u)-E_u(u^{-1}\dot
u)+T_u(u^{-1}u_x)+E_u(u^{-1}u_x)\Big)\ ,} and as we know
$s_xs^{-1}$ and $(T_u+E_u)(u^{-1}u_x)$, we can find
$(T_u-E_u)(u^{-1}u_t)$. From this we can in principle find
$u^{-1}\dot u$ as the function $\xi\mapsto T_u(\xi)-E_u(\xi)$ is
1-1 (if $\eta$ lay in the kernel of this operator then
$\eta+T_u(\eta)=
\eta+E_u(\eta)\in u^{-1}({\cal E}_+\cap {\cal E}_-)u=\{0\}$).

If we have a system with coordinates for configuration space $q_i$,
and Lagrangian $L(q_i,\dot q_i)$, then the canonical momenta are
$p_i=\partial L/\partial \dot q_i$, and we define a symplectic form on the
phase space
by
$\omega=\sum dp_{i}\wedge dq_{i}$. With a little thought,
it can be seen that this corresponds to the directional
derivative formula (where we have taken a Lagrangian
density $\cal L$)
\[
\omega(u,\dot u;a,b;c,d)\ =\ \int^{\pi}_{x=0}\Big(
{\cal L}''(u,\dot u;0,c;a,b)-{\cal L}''(u,\dot u;0,a;c,d)\Big)\ dx\ .
\]

If we write a change in $k$ as labelled by $y$ we get $k_y=u_ys+us_y$,
and likewise for $k_z=u_zs+us_z$. From the last section,
 we can write the
Lagrangian density for our system as
\[
4{\cal L}(u,\dot u)\ =\ \< E_u(u^{-1}\dot u-u^{-1}u_x),u^{-1}\dot
u+u^{-1}u_x\> \ ,
\]
so we can calculate a partial derivative
\begin{eqnarray*}
4{\cal L}'(u,\dot u;0,c) &=&
 \< E_u(u^{-1}c),u^{-1}\dot u+u^{-1}u_x\> +
\< E_u(u^{-1}\dot u-u^{-1}u_x),u^{-1}c\>  \cr
&=& \< E_u(u^{-1}\dot u)-T_u(u^{-1}\dot u)-E_u(u^{-1} u_x)-T_u(u^{-1} u_x),
u^{-1}c\>\ ,
\end{eqnarray*}
so
$
2{\cal L}'(u,\dot u;0,u_y)= -\<s_xs^{-1},u^{-1}u_y\>
$,
which results in
\begin{eqnarray*}
2{\cal L}''(u,\dot u;0,u_y;u_z,\dot u_z) &=&
-\<(s_xs^{-1})_z,u^{-1}u_y\>+\<s_xs^{-1},u^{-1}u_zu^{-1}u_y\> \cr
&=&
-\<(s_zs^{-1})_x,u^{-1}u_y\>+\<[s_xs^{-1},s_zs^{-1}],u^{-1}u_y\>
+\<s_xs^{-1},u^{-1}u_zu^{-1}u_y\>
\end{eqnarray*}
Now compare this with the standard 2-form on the
loop group of $D$. Consider
\begin{eqnarray*}
\<(k^{-1}k_y)_x,k^{-1}k_z\> &=&
\<(s^{-1}s_y)_x+[s^{-1}u^{-1}u_ys,s^{-1}s_x]+s^{-1}(u^{-1}u_y)_xs,
s^{-1}s_z+s^{-1}u^{-1}u_zs\> \cr   &=&
\< (s_ys^{-1})_x-[s_xs^{-1},s_ys^{-1}],u^{-1}u_z\>+\<
[s_xs^{-1},s_zs^{-1}],u^{-1}u_y\> \cr &&
+\<s_xs^{-1},[u^{-1}u_z,u^{-1}u_y]\>+\<s_zs^{-1},(u^{-1}u_y)_x\>\ .
\end{eqnarray*}
On integration we find
\[
\Big[\<s_zs^{-1},u^{-1}u_y\>\Big]_{x=0}^\pi\ =\
\int_{x=0}^\pi \Big(
\<(s_zs^{-1})_x,u^{-1}u_y\>+\<s_zs^{-1},(u^{-1}u_y)_x\>\Big)\ dx\ ,
\]
so we have the following symplectic form on the phase space:
\eqn{symp}{
2\omega(k;k_z,k_y)\ =\ \int^{\pi}_{x=0}
\<(k^{-1}k_y)_x,k^{-1}k_z\>\ dx\ -\
\Big[\<s_zs^{-1},u^{-1}u_y\>\Big]_{x=0}^\pi\ .}
Now we come to the complication, the fact that this form is
degenerate on
$C^\infty\big((0,\pi),D\big)$. If we take a change in
$k\in C^\infty\big((0,\pi),D\big)$ given by $k\phi$
for $\phi\in \gm$, then $\omega(k;k_z,k\phi)=0$ for all
$k_z$. To remedy this we could remove the null direction
by declaring that the phase space would actually be
$C^\infty\big((0,\pi),D\big)/M$. Equivalently
 we could consider the phase space
to consist of those $k=us\in C^\infty\big((0,\pi),D\big)$
for which $s(0)$ is the identity in $M$.

\subsection{The Hamiltonian density}
The Hamiltonian density generating the time evolution can be calculated by
\[
4{\cal H}\ =\ 4{\cal L}'(u,\dot u;0,\dot u)- 4{\cal L}(u,\dot u)\ ,
\]
and using our previous result we can write this as
\begin{eqnarray*}
4{\cal H} &=& -\<E_u(u^{-1}\dot u-u^{-1}u_x),u^{-1}u_x\>-
\<s_xs^{-1}+E_u(u^{-1}\dot u-u^{-1}u_x),u^{-1}\dot u\>  \cr
&=&
 -\<E_u(u^{-1}\dot u-u^{-1}u_x),u^{-1}u_x\>-
\<
T_u(u^{-1}\dot u)+T_u(u^{-1}u_x),u^{-1}\dot u\>  \cr
&=&
 \<E_u(u^{-1}u_x),u^{-1}u_x\>-
\<
T_u(u^{-1}\dot u),u^{-1}\dot u\>\ =\  \<E_u(u^{-1}u_x),u^{-1}u_x\>+
\<
E_u(u^{-1}\dot u),u^{-1}\dot u\>\ ,
\end{eqnarray*}
or equivalently
\eqn{uhamilt}{
8{\cal H} \ =\  \<(E_u-T_u)(u^{-1}u_x),u^{-1}u_x\>+
\<(E_u-T_u)(u^{-1}\dot u),u^{-1}\dot u\>\ .}
Using the equation we derived for $s_xs^{-1}$,
we can rewrite $\<(E_u-T_u)(u^{-1}\dot u),u^{-1}\dot u\>$ as
\begin{eqnarray*} &&
\< (T_u+E_u)(u^{-1}u_x)-2s_xs^{-1},(E_u-T_u)^{-1}\Big(
(T_u+E_u)(u^{-1}u_x)-2s_xs^{-1}\Big)\> \\
&=& -\< (T_u+E_u)(E_u-T_u)^{-1}(T_u+E_u)(u^{-1}u_x),u^{-1}u_x\>
\\  && -\ 4
\<s_xs^{-1},(E_u-T_u)^{-1}(T_u+E_u)(u^{-1}u_x)\>
\ +\ 4 \<s_xs^{-1},(E_u-T_u)^{-1}(s_xs^{-1})\>
\end{eqnarray*}
If we observe that
\[
\<(E_u-T_u)(u^{-1}u_x),u^{-1}u_x\>\ =\
\<(E_u-T_u)(E_u-T_u)^{-1}(E_u-T_u)(u^{-1}u_x),u^{-1}u_x\>
\]
then we can write
\begin{eqnarray}\label{ushamilt}
4{\cal H}&=&
 -\< T_u(E_u-T_u)^{-1}E_u(u^{-1}u_x),u^{-1}u_x\>
 -\< E_u(E_u-T_u)^{-1}T_u(u^{-1}u_x),u^{-1}u_x\> \nonumber\\  && - \ 2
\<s_xs^{-1},(E_u-T_u)^{-1}(T_u+E_u)(u^{-1}u_x)\>
\ +\ 2 \<s_xs^{-1},(E_u-T_u)^{-1}(s_xs^{-1})\>\ .
\end{eqnarray}

To simplify this equation we shall first look at the form of the
projections to the subspaces $u^{-1}{\cal E}_+u$ and $u^{-1}{\cal
E}_-u$ in terms of the graph coordinates. If we take $\xi\in \gg$
and $\phi\in\gm$, we can write
\[
\xi+\phi\ =\ (w+E_u(w))\ +\ (y+T_u(y))\ ,
\]
where $w=(E_u-T_u)^{-1}\phi-(E_u-T_u)^{-1}T_u(\xi)$ and
$y=(E_u-T_u)^{-1}E_u(\xi)- (E_u-T_u)^{-1}\phi$. Then we can define
projections $\pi_{u+}$ and $\pi_{u-}$ to $u^{-1}{\cal E}_+u$ and
$u^{-1}{\cal E}_-u$ as
\[
\pi_{u+}(\xi+\phi)\ =\ w+E_u(w)\quad{\rm and}\quad \pi_{u-}(\xi+\phi)\
=\ y+T_u(y)\ .
\]

It follows that
\eqn{projxi}{
(\pi_{u+}-\pi_{u-})\xi=-2E_u(E_u-T_u)^{-1}T_u\xi\
- (E_u-T_u)^{-1}(T_u+E_u)\xi,}
\eqn{projphi}{
(\pi_{u+}-\pi_{u-})\phi\ =\
2(E_u-T_u)^{-1}\phi+(T_u+E_u)(E_u-T_u)^{-1}\phi\ .}
{}From this we can rewrite the last equation for the Hamiltonian
as
\[
4{\cal H}\ =\ \<  (\pi_{u+}-\pi_{u-})(u^{-1}u_x+s_xs^{-1}),
u^{-1}u_x+s_xs^{-1}\>\ .
\]
This can be further simplified by removing the $u$ dependence from the
projections. If $\pi_+$ is the projection to ${\cal E}_+$
with kernel ${\cal E}_-$, then $\pi_{u+}={\rm Ad}_{u^{-1}}\circ \pi_+\circ
{\rm Ad}_{u}$, and since the inner product is adjoint invariant we find
\eqn{hamilt}{
4{\cal H}\ =\ \<  (\pi_{+}-\pi_{-})(u_xu^{-1}+us_xs^{-1}u^{-1}),
u_xu^{-1}+us_xs^{-1}u^{-1}\>}
or in terms of combined variable on $D$,
\eqn{khamilt}{4{\cal H}=\ \<  (\pi_{+}-\pi_{-})(k_xk^{-1}),
k_xk^{-1}\>\ .}
The equations of motion can similarly be written
in terms of $k$ as
\eqn{keqm}{\dot k k^{-1}=(\pi_--\pi_+)(k_x k^{-1}).}

\subsection{Symmetries of the models}
Returning to the equations of motion in the form
$k_\pm k^{-1}\in{\cal E}_\mp$, it
is clear that
\eqn{rightsym}{ k\mapsto kd,\quad d\in D}
is a global symmetry of the model. This has been discussed in
\cite{KliSev:poi}.
In addition to this known symmetry we now consider
\eqn{leftsym}{ k\mapsto dk,\quad{\cal E}_\pm\mapsto d{\cal E}_\mp d^{-1},
\quad d\in D}
which alters the subspaces ${\cal E}_\pm$ and hence the model. On
our phase space picture, where the different subspaces appear as
different Hamiltonians, this left translation in $D$ may not
preserve the Hamiltonian for a particular model, but rather takes
us from one model to another.

To have a dynamical symmetry of a particular model we can proceed
to restrict to left multiplication by those $d\in D$ such that $d{\cal
E}_\pm d^{-1}={\cal E}_\pm$. We distinguish two special cases: (1)
The subspaces ${\cal E}_\pm$ are $G$-invariant, and (2) The
subspaces ${\cal E}_\pm$ are $M$-invariant. In case (1) we say
that the models are {\em $G$-invariant}. Then $T_u=T_e$ and
$E_u=E_e$ are independent of $u\in G$, and the models themselves
are simpler to work with. The actions of $d\in G$ by left
translation in terms of the variables of the model and the dual
model are
\[ (u,s)\mapsto (du,s),\quad (t,v)\mapsto ((t^{-1}\ra d^{-1})^{-1},
(t^{-1}\la d^{-1})^{-1}v)\] respectively.
To see if the left translation has a moment map,
we consider $k_z=\delta k$ for $\delta\in\gd$ in the equation for the
symplectic form:
\[
2\omega(k;\delta k,k_y)\ =\ \int^{\pi}_{x=0}
\<k(k^{-1}k_y)_xk^{-1},\delta\>\ dx\ -\
\Big[\<s_zs^{-1},u^{-1}u_y\>\Big]_{x=0}^\pi\ .
\]
If $\delta\in \gg$, then $s_z=0$, so we have the moment map
\[ I_\delta(k)=-\h\int \<k_x k^{-1},\delta\> \extd x\ ,\quad \delta\in\gg.\]
In terms of the sigma-model on $G$, this is
\[ -4I_\delta(u)=\int \<2u^{-1}u_x+(T_u-E_u)(u^{-1}\dot u)
+(T_u+E_u)(u^{-1}u_x),u^{-1}\delta u\>\extd x\
,\quad\delta\in\gg\] which is a conserved charge in the
$G$-invariant case. The left translations for $\delta\in\gm$ are
not in general given by moment maps.

There are analagous formulae for the dual model and the $M$-invariant case.
We shall return to these symmetries when we have have discussed
boundary conditions for the models.
We shall also study the particular properties of $G$-invariant models in
some detail in later sections.

\section{Solutions independent of $x$}

In this section we show that the systems above in the Hamiltonian
form have `point-like' limits where the solutions are restricted
so that the field $u$, say, is independent of $x$. This then
becomes a system of a classical particle moving on the group
manifold of $G$. In the dual picture, i.e. in terms of the
variable $t$, the model is far from point-like and instead
describes some form of extended object in the manifold $M$. We
obtain the Poisson brackets and the Hamiltonian and we study the
symmetries, in particular the $G$-invariant case. The dual case
where $t$ is pointlike and $u$ extended is identical with the
roles of $G$ and $M$ interchanged and is therefore omitted except
with regard to the study of this case when the model is
$G$-invariant.

\subsection{The point-particle Poisson structure}

The solutions which have $u(x)$ independent of $x$
 are parameterised by
initial values of $u\in G$ and $p=s_xs^{-1}\in\gm$. This is
because the equation $s_xs^{-1}=(T_u-E_u)(u^{-1}\dot u)/2$ shows
that $p$ is also independent of $x$. Therefore the effective phase
space coordinates are $(u,p)$ rather than the fields $(u(x),s(x))$
in the general case. The symplectic form per unit length is then
\[
2\omega(u,p;u_z,p_z,u_y,p_y)\ =\ \<p_y,u^{-1}u_z\>-\<p_z,u^{-1}u_y\>+
\<p,[u^{-1}u_z,u^{-1}u_y]\>\ ,
\]
which is closed independently of the pairing used. This can also
be written as
\eqn{consym}{2\omega(u,p;u_z,p_z;u_y,p_y)
=\<(upu^{-1})_y,u_zu^{-1}\>-\<(upu^{-1})_z,u_yu^{-1}\>
-\<p,[u^{-1}u_z,u^{-1}u_y]\>.}

We now invert
the symplectic form on the phase space $\gm\times G$
 to find the Poisson structure.
Define $\omega_0:(\gm\oplus\gg)\tens (\gm\oplus\gg)\to\Bbb R$ by
\[
2\omega_0(p_y\oplus\xi_y,p_z\oplus\xi_z)\ =\
\<p_y,\xi_z\>-\<p_z,\xi_y\>+
\<p,[\xi_z,\xi_y]\>\ ,\quad\forall p_y,p_z\in\gm,\
\xi_y,\xi_z\in\gg
\]
Take a basis $e_i$ of $\gg$ and a dual basis $e^i$ of $\gm=\gg^*$
(for $1\le i\le n$). Then we can take a basis of $\gm\oplus\gg$ as
$f_i=e^i$ for $1\le i\le n$ and $f_i=e_{i-n}$ for $n+1\le i\le
2n$. Then in this basis,
\[
2\omega_0
\ =\ \left(\begin{array}{cc} 0 & \id \cr -\id & A \end{array}\right)\quad
{\rm and}\quad (2\omega_0)^{-1}
\ =\ \left(\begin{array}{cc} A & -\id \cr \id & 0 \end{array}\right)
\]
where $A_{ij}=\<p,[e_i,e_j]\>$. The corresponding tensor is
\[
\h\omega_0^{-1}=\sum_{1\le i\le n}\Big( e_i\tens e^i-e^i\tens
e_i\Big)\ +\ \sum_{1\le i,j\le n}
\<p,[e_i,e_j]\>e^i\tens e^j\ .
\]
Now, $\omega(u,p;u\xi_z,p_z;u\xi_y,p_y)=\omega_0(p_z\oplus
\xi_z,p_y\oplus \xi_y)$ so its inverse, the corresponding Poisson
bivector is given by left translation from $\omega_0^{-1}$,
\eqn{bivec}{
\gamma(p,u)=2 \sum_i  \tilde e_i\tens e^i-e^i\tens\tilde e_i+2\delta p}
where $\tilde\xi=u\xi$ is the left-invariant vector field
generated by $\xi\in \gg.$

The Poisson bracket itself then can be described simply for
functions $f,g$ on $G$ and $\xi,\eta\in\gg=\gm^*$ by
\eqn{conpoi}{\{f,g\}=0,\quad \{\xi,f\}=-2\tilde\xi(f),\quad
\{\xi,\eta\}=2[\xi,\eta].}
{}From this it is clear that we can quantise the system with the
Weyl algebra $\C[G]\lcross U(\gg)$ or at the $C^*$-algebra level
$C(G)\lcross C^*(G)$ where $G$ acts on $G$ by left multiplication.

\subsection{The point-particle Hamiltonian}
We have shown that $p=s_xs^{-1}$ is independent of $x$, so $s$ is of the form
 $s=e^{px}a$, where $a\in M$
is also independent of $x$.
To find the equations of motion
we write $k=ue^{px}a$, where $u\in G$
depends only on time, not on $x$. Then the equation of motion
$\dot kk^{-1}=(\pi_--\pi_+)k_xk^{-1}$ gives
\[
\dot uu^{-1}\,+\, u\frac{d}{dt}(e^{px})e^{-px}u^{-1}
\,+\, ue^{px}\dot aa^{-1}e^{-px}u^{-1}\ =\ (\pi_--\pi_+)(upu^{-1})\ ,
\]
which yields, for the case $x=0$,
\[
u^{-1}\dot u
\,+\, \dot aa^{-1}\ =\ (\pi_{u-}-\pi_{u+})p\ ,
\]
and taking the first order terms in $x$ gives
\[
\dot p\ =\ [\dot aa^{-1},p]\ .
\]
We can now get rid of the variable $a$ and write the equations of motion in
terms of
$u$ and $p$ only,
\[
u^{-1}\dot u
\ =\ \pi_\gg(\pi_{u-}-\pi_{u+})p\ ,\quad \dot p\ =\
[\pi_\gm(\pi_{u-}-\pi_{u+})p,p]\ .
\]
In the constant case, the Hamiltonian per unit length (\ref{hamilt}) restricts
to
\eqn{conhamilt}{
4{\cal H}\ =\ \< (\pi_+-\pi_-)(upu^{-1}),upu^{-1}\>\ .}
We have to check that the restricted Hamiltonian and the
restricted symplectic form indeed correspond to these equations of
motion, i.e. that the constraint of $x$-independence commutes with
the original Hamiltonian. To do this, it will be convenient to first calculate
from the equations of motion
\[
\frac{d}{dt}(upu^{-1})\ =\ u[(\pi_{u-}-\pi_{u+})p,p]u^{-1}\ =\
[(\pi_{-}-\pi_{+})(upu^{-1}),upu^{-1}]\ ,
\]
and now we can write
\align{&&\equad 2\omega(u,p;u_z,p_z;\dot u,\dot p)\\
&&=\<[upu^{-1},(\pi_+-\pi_-)upu^{-1}],u_zu^{-1}\>-\<(upu^{-1})_z,\dot
u u^{-1}\>-\<upu^{-1},[u_zu^{-1},\dot u u^{-1}]\>\\
&&=\<[u_zu^{-1},upu^{-1}],(\pi_+-\pi_-)upu^{-1}\>
-\<(upu^{-1})_z-[u_zu^{-1},upu^{-1}],\dot u u^{-1}\>\\
&&=\<(upu^{-1})_z,(\pi_+-\pi_-)upu^{-1}\>
-\<up_zu^{-1},(\pi_+-\pi_-)upu^{-1}\>-\<up_zu^{-1},\dot u
u^{-1}\>\\
&&=\<(upu^{-1})_z,(\pi_+-\pi_-)upu^{-1}\>
-\<p_z,u^{-1}\dot u-(\pi_{u-}-\pi_{u+})p
\>\\
&&=\<(upu^{-1})_z,(\pi_+-\pi_-)upu^{-1}\>
+\<p_z,\dot aa^{-1}\>\ =\ 2{\cal H}_z\ ,}
where we used at the end the equations of
motion again,
and then that $\<p_z,\dot aa^{-1}\>=0$ as $\gm$ is isotropic.
In terms of graph coordinates, we can write the equations of motion as
\eqn{coneqmu}{
u^{-1}\dot u
\ =\ -2(E_u-T_u)^{-1}p=2T_u^{-1}(E_u^{-1}-T_u^{-1})^{-1}E_u^{-1}p}
\eqn{coneqmp}{\dot p\ =\ -[(E_u+T_u)(E_u-T_u)^{-1}p,p]
=[(E_u^{-1}-T_u^{-1})^{-1}(E_u^{-1}+T_u^{-1})p,p]}
and the Hamiltonian as
\eqn{conhamiltgraph}{
4{\cal H}\ =\ \< (\pi_{u+}-\pi_{u-})p,p\>\ =\ 2\< (E_u-T_u)^{-1}p,p\>
=2\<(E_u^{-1}-T_u^{-1})^{-1}E_u^{-1}p,E_u^{-1}p\>\ .}

There is also a `conjugate' description of the system which
we mention briefly here. Although only $s_xs^{-1}=p$ is
directly needed for solving the $x$-independent equations of
motion for the $u$ variable, the rest of the degrees of freedom in
$s$ are also an auxiliary part of the system from the point of
view of the the group $D$. It turns out
that one could equally regard $(p,a)$ as  phase space variables
and solve the system in terms of them, with $u$ regarded as
auxiliary. Then the equations of motion would be
\eqn{paeqm}{ \dot a a^{-1}
\ =\ \pi_\gm(\pi_{u-}-\pi_{u+})p=-(E_u+T_u)(E_u-T_u)^{-1}p\ ,
\quad \dot p\ =\
[\pi_\gm(\pi_{u-}-\pi_{u+})p,p]\ .}
If we work with the phase space $\gm\times M=\cg^\star\otimes
G^\star$, we can more easily compare the system with the classical
phase space of the bicrossproduct Hopf algebra $U(\cg)\cobicross
\C[G^\star]$ associated to the same factorisation of $D$ in
\cite{Ma:pla}. In fact both the Poisson structures and the natural
Hamiltonians look somewhat different, but the general
interpretation as a particle on $M=G^\star$ with momentum given by
$p\in\cg^\star$ is the same.

\subsection{Symmetries of the point-particle system}

We now consider which of the translation symmetries of the general
theory restrict to the $x$-independent solutions. First of all,
the right translation symmetries are not interesting in this case:
the right action by $M$ is the identity on our $(u,p)$ coordinates,
 while the right action by $G$ does not preserve that
$u$ is $x$-independent.
On the other hand, the left translation symmetries by $d\in D$ do
preserve that $u$ is $x$-independent. We compute the Hamiltonian
functions for these actions. First of all, for an infinitesimal
transformation by $\phi\in \gm$ the variations of $u,upu^{-1}$ are
\[ u_\phi=\phi\la u,\quad (upu^{-1})_\phi=[\phi,upu^{-1}]\]
and hence (\ref{consym}) yields
\[ 2\omega(u,p;u_z,p_z;u_\phi,p_\phi)=-\<(upu^{-1})_z,\phi\>\]
for any variation $u_z,p_z$. Hence the Hamiltonian function generating
this flow is
\[ I_\phi(u,p)=-\h\<upu^{-1},\phi\>=-h\<u(p\la u^{-1}),\phi\>=-\h
\<ub_p(u^{-1})u^{-1},\phi\>,\quad\forall\phi\in\gm.\]

Similarly, for an infinitesimal left translation generated by $\xi\in\gg$
we have $u_\xi=\xi u$ (the right-invariant vector field generated by $\xi$)
and $p_\xi=0$. In this case we obtain more simply
\[ 2\omega(u,p;u_z,p_z;u_\xi,p_\xi)=-\<(upu^{-1})_z,\xi\>\]
or the generating function
\[ I_\xi(u,p)=-\h\<upu^{-1},\xi\>=-\h\<p\ra u^{-1},\xi\>,
\quad\forall\xi\in\gg.\]
The two cases can be combined into a single generating function
or moment map
\eqn{conmom}{ I_\delta(u,p)=-\h\<upu^{-1},\delta\>,\quad
\forall \delta\in\gd.}

In particular, we see that if the model is $G$-invariant, so that
$G$ is a dynamical symmetry, then the projection of $upu^{-1}$ to
$\gm$,
\eqn{coninvmom}{ Q_G=p\ra u^{-1}}
is a constant of motion, the conserved charge for the symmetry.
Likewise, if the model is $M$-invariant then the projection of
$upu^{-1}$ to $\gg$,
\eqn{concoinvmom}{ Q_M= ub_p(u^{-1})u^{-1}}
is a constant of motion.

The Hamiltonian and the equations of motion also simplify in the
$G$-invariant case, namely (\ref{coneqmu})-(\ref{conhamiltgraph})
with $E_u=E_e$ and $T_u=T_e$. Writing $U=2(T_e-E_e)^{-1}$,
$V=\h(E_e+T_e)$, we have
\eqn{coninveqm}{ u^{-1}\dot u=Up\ ,\quad \dot p=[VU p,p]\ , \quad
4{\cal H}=-\<  Up,p  \> \ .}
Thus, the equations of motion decouple in this case; $\dot p$ is a
quadratic function of $p$ and $u^{-1}\dot u$ is a linear function
of $p$, i.e. can then be obtained (in principle) by integrating
$p(t)$.

\subsection{The extended system dual to the point-particle limit}

The dual model when $u$ is $x$-independent is described by
variables $t,v$ both far from $x$-independent.  The dual
 constraint
is one where $t$ is fixed to be $x$-independent, in which case
 the model
in our original $u,s$ description is far from $x$-independent.
Rather,
it is some form of `extended solution'.

We can reverse the order of factorisation $k=ue^{px}a=tv$
to get $t^{-1}=(e^{px}a)^{-1}\ra u^{-1}$ and
$v^{-1}=(e^{px}a)^{-1}\la u^{-1}$. Here $u$, $p$ and $a$ are
functions of $t$ only. It can be seen that $t$ has a modified
exponential behaviour in $x$, and that $v$ is a constant
acted on by an exponential as a function of $x$. In particular
$t$ will not satisfy the Neumann boundary conditions.

 The Hamiltonian can be written as
\[
4{\cal H}\ =\
\< (\pi_{t+}-\pi_{t-})(t^{-1}t_x+v_xv^{-1}),t^{-1}t_x+v_xv^{-1}
\>\ ,
\]
where $\pi_{t\pm}$ are the projections to $t^{-1}{\cal E}_\pm t$.
The constraints on the dual system corresponding to
constant $u$ are that $t\la v$ and $t_xt^{-1}+tv_xv^{-1}t^{-1}$
are independent of $x$.

\section{More about graph coordinates}

In this section we provide some preliminary results on the explicit
construction of the graph coordinates of the subspaces
${\rm Ad}_{u^{-1}}{\cal
E}_\pm$ in terms of the actions of the groups on the Lie algebras. This is
needed, in particular, for the explicit computations for the
quasitriangular case in the next section.  In fact it will be convenient to
consider the inverses of the graph coordinates rather than the graph
coordinates
themselves, as the formulae are considerably simpler.

Thus, given generic $\CE_\pm$, the subspace ${\rm Ad}_{u^{-1}}{\cal E}_+$
contains elements of the form
\[
{\rm Ad}_{u^{-1}}(E_e^{-1}(\phi)\oplus\phi)\ =\
{\rm Ad}_{u^{-1}}(E_e^{-1}(\phi)+b_\phi(u))\oplus \phi\ra u\ =\
E_u^{-1}(\phi\ra u)\oplus \phi\ra u\ ,
\]
so we deduce that
\[
E_u^{-1}(\phi)\ =\ {\rm Ad}_{u^{-1}}(E_e^{-1}(\phi\ra
u^{-1})+b_{\phi\ra u^{-1}}(u)).
\]
We can write this as
\eqn{barEb}{\bar E_u^{-1}\equiv \Ad_u\circ E_u^{-1}\circ((\ )\ra u),\quad
\bar E_u^{-1}=E_e^{-1}+b(u).}  Also observe that $(s\ra u^{-1})\la u=(s\la
u^{-1})^{-1}$ for any double cross product group, which implies
that ${\rm Ad}_{u^{-1}}b_{\phi\ra u^{-1}}(u)=
-b_\phi(u^{-1})$ (this is part of the cocycle property for $b$).
Hence we can write equivalently
\eqn{Einvb}{
E_u^{-1}(\phi)\ =\ {\rm Ad}_{u^{-1}}(E_e^{-1}(\phi\ra u^{-1}))\ -\
b_\phi(u^{-1})\ .} The same formulae hold for $T$ replacing $E$.

If we consider the dual model the subspace ${\rm Ad}_{t^{-1}}{\cal
E}_+$ contains elements of the form
\[
{\rm Ad}_{t^{-1}}(\xi\oplus\hat E_e^{-1}(\xi))\ =\ t^{-1}\la \xi \oplus
{\rm Ad}_{t^{-1}}(\hat E_e^{-1}(\xi)+a_\xi(t^{-1}))\ =\ t^{-1}\la \xi \oplus
\hat E_t^{-1}(t^{-1}\la \xi)\ ,
\]
from which we deduce
\eqn{hatbarEa}{\hat{\bar E}^{-1}_t\equiv
\Ad_{t}\circ \hat{E}_t^{-1}\circ(t^{-1}\la(\ )),\quad \hat{\bar E}^{-1}_t
=\hat{E}_e^{-1}+a(t^{-1})}
or equivalently that
\eqn{hatEinva}{\hat E_t^{-1}(\xi)\ =\ {\rm Ad}_{t^{-1}}(\hat E_e^{-1}(t\la 
\xi))\ -\
a_\xi(t)\ .} Similarly for $\hat T$. Note also that
$E_e^{-1}(\phi)+\phi\in{\cal E}_+$ for all $\phi\in\gm$ and since
this also characterises $\hat E_e$ (and similarly for $\hat T_e$),
we conclude that
\eqn{EEinv}{ \hat E_e=E_e^{-1},\quad \hat T_e=T_e^{-1}.}

Finally, we specialise to the case of a coadjoint matched pair, i.e.
where $\gg$ is
a Lie bialgebra and $\gm=\gg^\star$, with $\gd=\gg\dcross\gg^\star$ the
Drinfeld double.
Now, associated to the Lie bialgebra structure is a Poisson-Lie group structure
on $G$ defined by bivector
\[ \gamma_G(u)= \tilde{\Pi(u)}\]
where $\tilde{\ }=R_*$ denotes extension as a left-invariant
vector field and $\Pi:G\to \gg\tens\gg$ is the cocycle $\Pi\in
Z^1_{\Ad}(G,\gg\tens\gg)$ extending the Lie cobracket
$\delta\in\Z^1_{\rm ad}(\gg,\gg\tens\gg)$ (which is the derivative
of $\Pi$ at the group identity). Since the action of $\gg^\star$
on $\gg$ in the coadjoint matched pair is just $\delta$ viewed by
evaluation against the second factor of its output, the cocycle
generator $b$ of its corresponding vector fields on $G$ us just
$b=\Pi$ in this case. Also observe that we could equally well have
defined $\gamma$ as generated by {\em right}-invariant vector
fields from some $\Pi^R$, say. Here
\[ \Pi^R(u)=\Ad_{u^{-1}}(\Pi(u))=-\Pi(u^{-1}),\]
the last equation by the cocycle condition obeyed by $\Pi$.

To apply these observations to the above we write operator
$E_u^{-1}:\gm\to\gg$ as evaluation against the {\em second} factor
of elements $E_u^{-1}\in\gg\tens\gg$ (we use the same symbols when
the meaning is clear). Similarly for $\hat E_t^{-1}$. Then
\eqn{EPi}{
E_u^{-1}=\Ad_{u^{-1}}(E_e^{-1})-\Pi(u^{-1})=\Ad_{u^{-1}}(E_e^{-1})+\Pi^R(u)}
as elements of $\gg\tens\gg$. Inverting this defines the
Lagrangian for our models,
\eqn{lagE}{ {\cal
L}=\<E_u(u^{-1}u_-),u^{-1}u_+\>=E_u(u^{-1}u_+,u^{-1}u_-)} where in
the second expression we view $E_u:\gg\to\gm$ as evaluation
against the second factor of $E_u\in\gm\tens\gm$. Or in terms of
$\bar E_u^{-1}=\Ad_u(E_u^{-1})\in\gg\tens\gg$, we have
\eqn{barEpi}{ \bar E_u^{-1}=E_e^{-1}+\Pi(u)}
and the Lagrangian written equally as
\eqn{lagbarE}{ {\cal L}=\<\bar E_u(u_-u^{-1}),u_+u^{-1}\>=\bar
E_u(u_+u^{-1},u_-u^{-1}).} One or other of these two forms is
usually easier to compute.

Similarly, for the dual model we identify $a(t):\gg\to\gm$ with
evaluation against the first component of $\hat\Pi^R$, i.e.
$a=-\hat
\Pi^R$ when the latter is considered as an operator by evaluation
against its second factor (a convention that we adopt unless
stated otherwise). Then
\eqn{hatEpi}{  \hat E_t^{-1}=\Ad_{t^{-1}}(E_e)+\Pi^R(t),\quad \hat {\bar
E}_t^{-1}=E_e+\hat\Pi(t)} and
\eqn{hatLagE}{ {\cal L}=\hat E_t(t^{-1}t_+,t^{-1}t_-)=\hat{\bar
E}_t(t_+t^{-1},t_-t^{-1})} is the Lagrangian for the dual model.

These results allow us to explicitly construct the graph
coordinates and the Lagrangians given a generic splitting of $\gd$
into subspaces $\CE_\pm$. The latter are equivalent to specifying
$E_e^{-1},T_e^{-1}$ and these allow us to obtain the general
$E_u^{-1}$ etc., from (\ref{Einvb}) or from (\ref{EPi}) etc., in the
coadjoint case.

\section{Models based on $\gg$ quasitriangular}

In this section we define a class of Poisson-Lie dual models based
on the double of $\gg$ (the usual setting) but in the special case
where $\gg$ is quasitriangular and factorisable. In this case were
are able to obtain much more explicit formulae for the model and
the dual model than in the general case.

A
Lie bialgebra is quasitriangular if there is an element
$r\in\gg\tens\gg$ such that $\delta\xi=\ad_\xi(r)$ and $r$ obeys
the classical Yang-Baxter equations
\eqn{cybe}{[r_{12},r_{13}]+[r_{12},r_{23}]+[r_{13},r_{23}]=0}
and has $2r_+=r+r_{21}$ ad-invariant. A factorisable
quasitriangular Lie bialgebra is one where $2r_+$ viewed as a map
$\gg^*\to\gg$ is invertible. We denote its inverse by $K$. In
standard examples where $\gg$ is simple, $K$ is a multiple of
the Killing form viewed as a map.

In
this case there is an isomorphism \cite{Sem:wha}\cite{Ma:book}
\[ \gd=\gg\dcross\gg^*\isom \gg_L\codcross\gg_R,\quad
\xi\oplus\phi\mapsto (\xi+r_1(\phi),\xi-r_2(\phi))\]
which also sends the bilinear form $\<\ ,\ \>$ on $\gd$ to $K_L-K_R$ on the two
copies $\gg_L,\gg_R$ of $\gg$. Here $K_L,K_R$ are two copies of
$K$. Therefore the inverse image of $\gg_L,\gg_R$ define a splitting of
$\gd$ into mutually orthogonal
subspaces. From the explicit form of the isomorphism in \cite{Ma:book} 
one finds
\eqn{codsplit}{ \CE_+=\{\xi-r_1(K(\xi))+K(\xi)\},\quad
\CE_-=\{\xi-r_2(K(\xi))-K(\xi)\}.}
These subspaces are not generic, however (the graphs blow up) but
they are the model for the construction which follows. In fact one
has a two parameter family of models by varying the
coefficients of $r_1,K$ in $\CE_+$, etc., with graph coordinates in the
general case. In another degenerate limit of these parameters one has the
principal sigma model as well.

\subsection{Construction of the quasitriangular models on $G$}

The subspaces $\CE_\pm$ defining our model will be constructed
 by introducing parameters into  (\ref{codsplit}) in such a way as to
 preserve orthogonality.
Equivalently, one may define suitable $E_e^{-1},T_e^{-1}$. We then
obtain the general graph coordinates by the method of Section~5. In fact
we consider the second problem first as it leads to the most elegant
choice of ansatz for the $E_e^{-1}$ etc.

Thus, in the case of a quasitriangular Lie
bialgebra one has simply
\eqn{Pir}{\Pi(u)={\rm Ad}_u(r)-r} for the
cocycle defining its Poisson structure. This defines the Drinfeld-Sklyanin
bracket on $G$ when $\gg$ is the standard quasitriangular structure\cite{Dri}
for a simple Lie algebra $\gg$. These are also the Poisson brackets of which
the
associated quantum groups in this case are the quantisations. We refer to
\cite{Ma:book} for further discussion of these preliminaries. In view
of (\ref{Pir}) and the results of Section~5, it is then immediate that the
graph coordinates for the model on $G$ in the quasitriangular case obey
\eqn{Fgenqua}{
E^{-1}_u=\ {\rm Ad}_{u^{-1}}(E^{-1}_e-r)+r\ }
as an element of $\gg\tens\gg$. This equation, together with a little linear
algebra, allows the explicit
computation
of the graph coordinates for any model based on a quasitriangular Lie
bialgebra,
given suitable $F_e$.

Motivated by (\ref{codsplit}) we now let
\[ E^{-1}_e=(\lambda+1)r+\mu K^{-1}\]
where $\lambda,\mu$ are two complex parameters. For generic values we will
indeed be able to invert to obtain graph coordinates $E_u,T_u$ and hence
will obtain a model of the type studied in Sections~2,3.

Clearly, from (\ref{Fqua}), we have
\eqn{Fqua}{ E^{-1}_u=\lambda\Ad_{u^{-1}}(r)+r+\mu K^{-1}}
as solving the equation (\ref{Fgenqua}) for all $\lambda,\mu$. If we denote
by $r_2:\gg^*\to\gg$ the evaluation against the second factor of
$r\in\gg\tens\gg$
and similarly by $r_1$ for evaluation against the first factor, we have
equivalently, as maps $\gm\to\gg$,
\eqn{Einvqua}{ E_e^{-1}=(\lambda+1)r_2+\mu K^{-1}=(\lambda+\mu+1)r_2+\mu r_1}
for our class of models. Similarly,
\eqn{Tinvqua}{ T_e^{-1}=-(\lambda+1)r_1-\mu K^{-1}
=-(\lambda+\mu+1)r_1-\mu r_2.}
These imply
\eqn{ETqua}{E_e^{-1}-T_e^{-1}=(\lambda+1+2\mu)K^{-1},\quad
E_e^{-1}+T_e^{-1}=(\lambda+1)(r_2-r_1).}

For further computations in the Hamiltonian formulation we need the difference
of the associated
projectors $\pi_\pm$. Rearranging (\ref{projxi})--(\ref{projphi}), we have
\eqn{pixi}{ (\pi_{u+}-\pi_{u-})\xi=2(E^{-1}_u-T^{-1}_u)^{-1}\xi
+(E_u^{-1}+T_u^{-1})
(E_u^{-1}-T_u^{-1})^{-1}\xi,\quad\forall\xi\in\gg,}
\eqn{piphi}{(\pi_{u+}-\pi_{u-})\phi
=-2E_u^{-1}(E_u^{-1}-T_u^{-1})^{-1}T^{-1}_u\phi
-(E_u^{-1}-T_u^{-1})^{-1}(E_u^{-1}+T_u^{-1})\phi,\quad\forall\phi\in\gm.}
Evaluating at the identity and inserting the above results for $E^{-1}_e$,
etc.,
we obtain:
\eqn{hamqua1}{\<(\pi_+-\pi_-)\xi,\xi\>={2\over \lambda+1+2\mu}K(\xi,\xi),\quad
\<(\pi_+-\pi_-)\xi,\phi\>={\lambda+1\over \lambda+1+2\mu}K(\xi,(r_1-r_2)\phi)}
\eqn{hamqua2}{\<(\pi_+-\pi_-)\phi,\phi\>
={2\over \lambda+1+2\mu}K(T_e^{-1}\phi,T_e^{-1}\phi)}
\eqn{hamqua3}{
K(T_e^{-1}\phi,T_e^{-1}\phi)={(\lambda+1)^2\over 4}K((r_1-r_2)\phi,(r_1-r_2)
\phi)
+{(\lambda+1+2\mu)^2\over 4}K^{-1}(\phi,\phi).}
These results provide for the computation of Hamiltonian
from (\ref{hamilt}) in Section~3.

It remains to show that the above $E_e^{-1},T_e^{-1}$ indeed define
an orthogonal splitting of $\gd$ into subspaces ${\cal E}_\pm$ and
to give these explicitly. First of all the corresponding subspaces
defined by our choice of $E_e^{-1},T_e^{-1}$ are
\eqn{Epqua}{
{\cal E}_+=\{E_e^{-1}\phi\oplus\phi\}
=\{\xi-{(\lambda+1)r_1(K(\xi))-K(\xi)\over\lambda+1+\mu}:\xi\in\gg\},}
\eqn{Emqua}{
{\cal E}_-=\{T_e^{-1}\phi\oplus\phi\}
=\{\xi-{(\lambda+1)r_2(K(\xi))+K(\xi)\over\lambda+1+\mu}:\xi\in\gg\}.}
To show that these form an orthogonal decomposition of $\gd$, we calculate
the inner products
\[\<E_e^{-1}\phi\oplus\phi,T_e^{-1}\phi\oplus\phi\>
=\<E_e^{-1}\phi,\phi\>+\<\phi,T_e^{-1}\phi\>
=(\lambda+1)\<(r_2-r_1)(\phi),\phi\>=0,\]
\[ \<E_e^{-1}\phi\oplus\phi,E_e^{-1}\phi\oplus\phi\>
=\<E_e^{-1}\phi,\phi\>+\<\phi,E_e^{-1}\phi\>
=(\lambda+1+2\mu)K^{-1}(\phi,\phi),\]
\[ \<T_e^{-1}\phi\oplus\phi,T_e^{-1}\phi\oplus\phi\>
=\<T_e^{-1}\phi,\phi\>+\<\phi,T_e^{-1}\phi\>
=-(\lambda+1+2\mu)K^{-1}(\phi,\phi).\]
In particular, $\CE_\pm$ are mutually orthogonal as required (the latter
two equations show further that the inner product is nondegenerate on
each subspace). To show that the subspaces span
$\gd$ we need to show that
\[ \xi\oplus\phi=E^{-1}_e(\psi)+\psi+T_e^{-1}(\chi)+\chi\]
has a (unique) solution for $\psi,\chi\in\gm$ for all $\xi\in\gg$ and
$\phi\in\gm$.
Clearly $\psi+\chi=\phi$. Meanwhile, putting in the form of
$E_e^{-1},T_e^{-1}$ we
have
\[ \xi=\mu K^{-1}(\psi-\chi)+(\lambda+1)(r_2(\psi)-r_1(\chi))\]
which can be rearranged as
\[ \xi+(\lambda+1)(-r_2+\h K^{-1})(\phi)=\h(\lambda+1+2\mu)K^{-1}(\psi-\chi).\]
Thus we have an orthogonal splitting if and only if
\eqn{splitcond}{\lambda+1+2\mu\ne 0.}
We assume this throughout. Moreover, the splitting has the inverse-graph
coordinates $E_e^{-1},T_e^{-1}$ computed above.

This completes the construction of our model at least
in the Hamiltonian formulation. Indeed, this can be defined entirely
in terms of $E_u^{-1},T_u^{-1}$ without recourse to $E_u,T_u$ themselves. It is
clear from our construction that:

(1) The model is $G$-invariant if and only if
\eqn{invcond}{ \lambda=0}
(or the Lie bialgebra structure on $\gg$ is identically zero.)

(2) The standard Lagrangian for the model (which requires $E_u)$ exists if
and only if (\ref{Fqua}) are nondegenerate, in particular when $\mu
K$ dominates, i.e.
\eqn{lagrcond}{ |\mu|>>|\lambda+1|}
and $\gg$ is semisimple.

We describe several special cases.

\subsubsection*{Modified principal sigma model.} This is obtained by
$\lambda=-1,\mu=1$.
Then
\eqn{modpri}{ \CE_\pm=\{\xi\pm K(\xi):\ \xi\in\gg\},\quad E_e^{-1}
=K^{-1}=-T_e^{-1}}
Here $E_u$ is obtained by inverting $F_u=K^{-1}-\Pi(u)$ and is not independent
of $u\in G$. Considering $K,\Pi$ as maps $K,\Pi_2$ by evaluation against the
second component, we have
\[ E_u^{-1}-T_u^{-1}=2K^{-1},\quad E_u^{-1}+T_u^{-1}=2\Pi^R(u)\]
for this model. Here $\Pi^R(u)$ defines the Poisson-bracket associated
to the Lie bialgebra structure of $G$ and is viewed as a map $\gm\to\gg$ by 
evaluation
(as usual) against its second factor. In particular, the Lagrangian
is
\eqn{modprilagr}{ {\cal L}=\<(K^{-1}+\Pi^R(u))^{-1}u^{-1}u_-,u^{-1}u_+\>=
\<(K^{-1}+\Pi(u))^{-1}(u_-u^{-1}),u_+u^{-1}\>.}

This recovers the setting of \cite{KliSev:poi}, for example, as a special
case of our class of models. Note that the formulae for general $\mu$ but
$\lambda=-1$ are strictly similar, with $E_u=(\mu K^{-1}+\Pi^R(u))^{-1}$
in the Lagrangian instead.

\subsubsection*{Pure-quasitriangular and principal sigma model.}

The $G$-invariant models are obtained by
$\lambda=0$, $\mu=0$. In this case
\[ E_u^{-1}=E_e^{-1}=r_2+\mu K^{-1},\quad T_u^{-1}=T_e^{-1}=-r_1-\mu K^{-1}.\]
For the equations of motion we can use the equations
$u^{-1}u_-=E_e^{-1}(s_-s^{-1})$ and
$u^{-1}u_+=T_e^{-1}(s_+s^{-1})$ since the operators
$E_e^{-1}$
and $T_e^{-1}$ are defined as above, even though $E_e$ and $T_e$ may not be.
 Then the equations of motion are most conveniently
described as a sigma model for $s$, with equation
\[
(T_e^{-1}(s_+s^{-1}))_-\,-\,
(E_e^{-1}(s_-s^{-1}))_+\ =\ -[ E_e^{-1}(s_-s^{-1}),
T_e^{-1}(s_+s^{-1})]\ .
\]
We see that this case contains another sigma model on the dual group
which makes sense in the G-invariant case. Indeed, in the general $G$-invariant
case the variable $s$ may be considered to
have a complex parameter $\mu$, which makes this look very much like inverse
scattering
for the sigma model. Moreover, for generic $\mu$, the operators $E_e$ and
$T_e$ do exist,
and both $u$ and $s$ are described by sigma models.

The pure-quasitriangular model is the special case with $\mu=0$ as well.
In this case the subspaces $\CE_\pm$ are the ones in (\ref{codsplit}) 
corresponding to
the Drinfeld double as $\gg\codcross\gg$. This new class of models has
Hamiltonian
defined by
\[ \h\<(\pi_+-\pi_-)(\xi\oplus\phi),\xi\oplus\phi\>
=K(\xi,\xi)+K(\xi,(r_1-r_2)\phi)
+K(r_1(\phi),r_1(\phi)).\]

The principal sigma model is the limit with $\mu\to\infty$
and a suitable rescaling. It is on the boundary of
our moduli space of quasitriangular models. Then
\eqn{pri}{ \CE_+=\{\xi+\mu^{-1}(\xi-r_1\circ K(\xi))+K(\xi))\},\quad \CE_-=
\{\xi+\mu^{-1}(\xi-r_2\circ K(\xi))-K(\xi))\}}
and
\[ E_e=\mu^{-1}(K^{-1}+\mu^{-1}r_2)^{-1}
=\mu^{-1}K(1-\mu^{-1}r_2\circ K+\cdots)\]
Hence the Lagrangian is
\eqn{prilagr}{ \CL(u)=\<(\mu K^{-1}+r_2)^{-1}(u^{-1}u_-),u^{-1}u_+\>
=\mu^{-1}K(u^{-1}u_-,u^{-1}u_+)+\mu^{-2}K(r_2\circ
K(u^{-1}u_-),u^{-1}u_+)+\cdots}
which after an infinite renormalisation has leading term the usual principal
sigma model.

The equation of motion, to lowest order in $\mu^{-1}$, is
\[
K((u^{-1}u_+)_-\,+\, (u^{-1}u_-)_+)\ =\
\mu^{-1}(K(r_1K(u^{-1}u_+)_-\,+\, r_2K(u^{-1}u_1)_+)\ -\
[K(u^{-1}u_-),K(u^{-1}u_+)])\ +\cdots
\]
This is the usual principal sigma model equations of motion to lowest order 
in $\mu^{-1}$, namely
\[ (u^{-1}u_+)_-\,+\, (u^{-1}u_-)_+)=0.\]

\subsection{Quasitriangular models on $SU_2$.}

We now compute these models for the group $G=SU_2$ and for its other real form
$G=SL_2(\R)$. Actually, only the second of these is strictly real and
quasitriangular.
Thus, with a basis $\{H,X_\pm\}$ for its Lie algebra (with the usual
relations), we
take the Drinfeld-Sklyanin quasitriangular structure
\[ r=X_+\tens X_-+{1\over 4}H\tens H.\]
Let $sl_2(\R)^\star$ have the dual basis $\{\phi,\psi_\pm\}$ then its Lie
algebra structure is
\[ {}[\phi,\psi_\pm]=\h \psi_\pm,\quad [\psi_+,\psi_-]=0\]
and the other required maps are
\[ r_2\pmatrix{\phi\cr\psi_+\cr\psi_-}=\pmatrix{{1\over 4}H\cr0\cr X_+},
\quad K\pmatrix{H\cr X_+\cr X_-}=\pmatrix{2\phi\cr\psi_-\cr\psi_+}.\]

Note that if we take a different real form
\[ e_1={-\imath\over 2}(X_++X_-),\quad e_2={-1\over 2}(X_+-X_-),\quad
e_3={-\imath\over 2}H\]
then $[e_i,e_j]=\eps_{ijk}e_k$ (the real form $su_2$) but
\[ r=-\sum_i e_i\tens e_i+\imath (e_1\tens e_2-e_2\tens e_1)\]
is not real in this basis. If $\{f_i\}$ is a dual basis then
\[ r_2(f_j)=-e_j+\imath e_i\eps_{ij3},\quad K=-\h\id.\]
This means that although we can arrange for a completely real
Lie bialgebra $su_2$ in this basis (here the
Lie coalgebra is purely imaginary but we can rescale $r$ to make it
real) it is not a quasitriangular one over $\R$; the required $r$ if we
want to obey (\ref{cybe}) lives in the complexification. In the above 
conventions the
Lie algebra $sl_2^\star$ in the dual basis is imaginary,
\[{} [f_i,f_j]=\imath(\delta_{ik}\delta_{j3}-\delta_{jk}\delta_{i3})f_k.\]
The choice of basis $e_i^\star=-\imath f_i$ is its real form $su_2^\star$.

\subsubsection*{Modified principal sigma model on $SU_2$.}

To construct the model we will need $\Pi(u)=\Ad_u(r_-)-r_-$ quite
explicitly, where $r_-=\imath e_1\wedge e_2$ is the
antisymmetric part of $r$.  For our purposes we write $SU_2$ as elements
\[ u=\pmatrix{a&b\cr -\bar b&\bar a},\quad |a|^2+|b|^2=1.\]
Then working with the matrix representation $e_i={-\imath\over 2}\sigma_i$
given by the Pauli matrices it is easy to find
\[\Ad_{u^{-1}}(e_1)=\Re(a^2-b^2) e_1+\Im(a^2+b^2)e_2-2\Re(ab)e_3\]
\[\Ad_{u^{-1}}(e_2)=-\Im(a^2-b^2)e_1+\Re(a^2+b^2)e_2+2\Im(ab)e_3\]
and hence
\eqn{suPiR}{ \Pi^R(u)=2\imath e_1\wedge e_2|b|^2-e_3
\wedge e_1(a\bar b-\bar a b)-\imath e_2\wedge e_3(a\bar b+\bar a b)}
after a short computation, which is purely imaginary (as expected). 
Evaluating against
the second factor and regarding
as a matrix we have
\[  E_u^{-1}=K^{-1}+\Pi^R(u)=-2\pmatrix{1&-\imath |b|^2 &-\imath\Im(a\bar b)\cr
\imath|b|^2 & 1&
\imath\Re(a\bar b)\cr
\imath\Im(a\bar b)& -\imath\Re(a\bar b)&1}.\]
Here $E_u^{-1}(f_j)=E^{-1}_{ij}e_i$, where $(E^{-1}_{ij})$ is the matrix
shown. Note that we can write
\[ E^{-1}_{ij}=-2(\delta_{ij}+\imath\eps_{ijk}\pi_k),\quad  \pi
=\pmatrix{\Re(a\bar b)\cr
\Im(a\bar b)\cr
-|b|^2 }\]
and any matrix of this form has inverse
\[ E_{ij}=-{1\over 2(1-\pi^2)}(\delta_{ij}-\imath\eps_{ijk}\pi_k-\pi_i\pi_j).\]
Here $\pi^2=\pi\cdot\pi=|b|^2$ in our case. The corresponding operator is
$E_u(e_j)=E_{ij}f_i$.
To cast the resulting
Lagrangian in a useful form let us note that
\[ \trace(\id-\pisl)\sigma_i\sigma_j=\trace(\id-\pi\cdot\sigma)(\delta_{ij}\id
+\imath\eps_{ijk}\sigma_k)=2(\delta_{ij}-\imath\eps_{ijk}\pi_k)\]
where $\sigma_i$ are the Pauli matrices and
$\pisl=\pi\cdot\sigma$. Hence in our representation of $su_2$ in
basis $e_i={-\imath\sigma_i\over 2}$ we have
\eqn{modprisu2lag}{ {\cal L}={1\over |a|^2}\left(\trace[(\id-\pisl)
u^{-1}u_+u^{-1}u_-]
-\h \trace[\pisl u^{-1}u_+]
\trace[\pisl u^{-1}u_-]\right)}
where
\[ \pisl=\pmatrix{-|b|^2&\bar a b\cr a\bar b&|b|^2}=\pmatrix{0&b\cr\bar b&0}u
=u^{-1}\pmatrix{0&b\cr\bar b&0}.\]

The matrix $E_{ij}$ here is complex since $\Pi^R$ in our conventions is
imaginary. For a completely real version of this model
on $SU_2$ one should keep the freedom of general $\mu$ in this class of
models so that $E_u^{-1}=\mu K^{-1}+\Pi^R(u)$ and then set $\mu=\imath$.
Taking the real normalisation of $su_2$ as a Lie bialgebra (i.e. multiplying
$r$ by $-\imath$ so that $r_-=e_1\wedge e_2$ and $K^{-1}=2\imath\id$) gives
the same $E_{ij}^{-1}$ as above but times $-\imath$ off the diagonal.
One may also work of course
on $G=SL_2(\R)$ with real $r,K$ for a completely real model with $\mu=1$.

This class of models has been considered specifically for $SU_2$ in 
\cite{LleVar:poi}, although not so explicitly as above.

\subsubsection*{Pure-quasitriangular and principal sigma models on $SU_2$.}

Here we take $\lambda=0$ and can write
down immediately
\[ E_u^{-1}=E_e^{-1}=-\pmatrix{1+2\mu&-\imath&0\cr \imath & 1+2\mu&0\cr
0&0&1+2\mu}\]
which has inverse
\[ E_u=E_e={-1\over 4\mu}\pmatrix{{1+2\mu\over1+\mu}&{i\over 1+\mu}&0\cr
{-\imath\over 1+\mu}&{1+2\mu\over1+\mu}&0\cr 0&0&{4\mu\over
1+2\mu}}\] for $\mu\ne 0,-\h,-1$. The Lagrangian defined by this
can be conveniently obtained by writing
\[ E_{ij}^{-1}=-(1+2\mu)(\delta_{ij}-\imath\eps_{ijk}\pi_k),\quad
\pi=\pmatrix{0\cr0\cr{1\over 1+2\mu}}\]
which implies (by similar computations to those above),
\[ {\cal L}={1\over\mu(1+\mu)}
\left(\trace[\pmatrix{1+\mu&0\cr 0&\mu}u^{-1}u_+u^{-1}u_-]
-{1\over 4(1+2\mu)}\trace[\sigma_3u^{-1}u_+]\trace[\sigma_3u^{-1}u_-]\right).\]
This singular for the pure quasitriangular model where $\mu=0$, and also does
not have a good limit at $\mu=\infty$ for the principal sigma model. Rather,
wehave well-defined equations of motion conveniently described as a sigma
model
for $s\in M$ as explained above, using $E_e^{-1}$ and a similar matrix for
$T_e^{-1}$.

On the other hand, by changing the normalisation of the Lie bialgebra structure
(namely, dividing $r$ by $\mu$) we have $E_u$  with the same matrix as above
but without the $\mu^{-1}$ factor in front. This rescaled Lagrangian
is well defined both for $\mu=0$ and $\mu=\infty$, with
\[ \mu {\cal L}\to
\cases{ \h\trace[(1+\sigma_3)u^{-1}u_+u^{-1}u_-]
-{1\over 4}\trace[\sigma_3u^{-1}u_+]\trace[\sigma_3u^{-1}u_-]&as\
$\mu\to 0$\cr
\trace[u^{-1}u_+u^{-1}u_-]&as\ $\mu\to\infty$}.\]
The first limit is the Lagrangian for the rescaled
pure-quasitriangular model on $SU_2$, while the second is the
standard Lagrangian for the  principal sigma model on $SU_2$ based
on the Killing form of $su_2$.

Notice that in this rescaled model the Lie cobracket of $su_2$
is infinite at $\mu=0$, i.e. the Lie algebra $\gm$ has infinite commutators,
and zero at $\mu=0$, i.e. the Lie algebra $\gm$ is Abelian. The geometrical
pictures behind these two models are therefore very different but interpolated
by general $\mu$.

Also note that the $\mu=0$ limit here is again defined by a complex
Lagrangian. For a real version one may look at the pure-quasitriangular model
on $G=SL_2(\R)$ instead. Here we have, clearly,
\[ E_e^{-1}\pmatrix{\phi\cr\psi_+\cr\psi_-}=\pmatrix{{1+2\mu\over 4}H\cr
\mu X_-\cr (1+\mu)X_+},\quad E_e\pmatrix{H\cr X_+\cr X_-}=\mu^{-1}
\pmatrix{{4\mu\over 1+2\mu}\phi\cr {\mu\over 1+\mu}\psi_-\cr\psi_+}.\]
As before, we take out a factor $\mu$ by rescaling in order to obtain
 well-defined operators
$E_e$ at $\mu=0,\infty$, this time with all coefficients being real in our
choice
of bases. The corresponding Lagrangian
can easily be written out explicitly upon fixing a description of
$u\in SL_2(\R)$. For example, if we write
\[ u=e^{xX_+}e^{hH}e^{yX_-}\]
so that
\[ u^{-1}u_\pm=x_\pm X_+e^{-2h}+(h_\pm + yx_\pm e^{-2h})H+(y_\pm-2y h_\pm
-2y^2x_\pm e^{-2h})X_-\]
using the relations of $sl_2$ then the rescaled $\mu=0$ limit gives the
Lagrangian
\[ {\cal L}=e^{-2h}x_+(y_--2yh_--2y^2e^{-2h}x_-)\]
as the pure-quasitriangular model on $SL_2(\R)$. The $\mu=\infty$ limit is
the standard principal sigma model on $SL_2(\R)$
and the general case interpolates the two.

\subsection{Dual of the quasitriangular models on $G^\star$}

The quasitriangular models are examples of the case where the
factorisation is based on the Drinfeld double associated to a Lie
bialgebra, so that $E_u^{-1}$ is related to the Poisson-Lie group
$G$. Hence the dual models are of the same form but based on the
Poisson-Lie group $G^\star$ rather than $G$, i.e. with with
$\hat\Pi(t)\in\gm\tens\gm$ in place of $\Pi$. As explained in
Section~5 we can then construct them from the initial data
\[ \hat E_e^{-1}=E_e,\quad \hat T_e^{-1}=T_e\]
as given above for our quasitriangular models. We compute
$\hat{\bar E}^{-1}_t=E_e+\hat\Pi(t)$ and invert to obtain the
Lagrangian
\eqn{dualbialag}{\hat{\cal L}=\<(E_e+\hat\Pi(t))^{-1}t_-t^{-1},t_+t^{-1}\>}
for the dual model. For the models below, where there is no
special $\Ad_t$-invariance of $E_e$, this is easier than computing
the Lagrangian via $\hat E_t$.

We outline the results for $SU_2^\star$ and $SL_2(\R)^\star$.
First of all we describe these groups explicitly. The former is
generated by the basis $\{-\imath f_i\}$, i.e. we write
$\phi=\phi_i(-\imath f_i)\in\gm$ for real $\phi_i$, which we
regard as a vector $\vec\phi$. One standard representation of the
resulting group is as matrices of the form
\[ \pmatrix{x&z\cr0 & x^{-1}},\quad x>0,\quad z\in\C.\]
This is the group occuring in the Iwasawa decomposition
$SL_2(\C)=SU_2\dcross SU_2^\star$,
see \cite{Ma:book}. Another description useful for very explicit computations
is as the semidirect product $\R^2\lcross\R$ \cite{Ma:book}, which can be
viewed
as a modified product on $\R^3$. Elements are $\vec s\in\R^3$ with $s_3>-1$
and the product law and inversion are
\[ \vec s\vec t=\vec s+(s_3+1)\vec t,\quad \vec s^{-1}=-{\vec s\over s_3+1}.\]
The exponentiation from the Lie algebra to a group is explicitly
\[ \vec s=\vec \phi {e^{\phi_3}-1\over\phi_3}\]
for $\vec s=e^{\vec\phi}$ in the natural 3-dimensional coadjoint
representation.
See \cite{Ma:book}. The real form $SL_2(\R)^\star$ has a similar description
as $\C\lcross\R$, i.e. where $s_2$ is imaginary and $s_1,s_3$ real
with $s_3>-1$
according to the conventions in \cite{Ma:book}. Note that $x=s_3+1$ is
multiplicative under the group law if one wants a more standard notation.

The Lie bracket on $su_2$ determines the Lie cobracket and Poisson
structure on $SU_2^\star$ (and similarly on $SL_2(\R)^\star$). It
is given by \cite{Ma:book}
\[ \hat\Pi(s)=-\imath(\eps_{ija}s_a+\h s^2\eps_{ij3})f_i\tens f_j.\]
Explicitly,
\[ \imath\hat\Pi(s)=\h(s_1^2+s_2^2+(s_3+1)^2-1)f_1\wedge f_2+s_2 f_3\wedge
f_1+s_1 f_2\wedge f_3.\]

Note also that the notation $s_\pm s^{-1}$ means more precisely
$R_{s^{-1}*}s_\pm$. Similarly for $s^{-1}s_\pm$. In our present
group coordinates, from the product law, it is easy to see that
\[ L_{s*}\vec\phi=(s_3+1)\vec\phi,\quad
R_{s*}\vec\phi=\vec\phi+\phi_3\vec s.\]

\subsubsection*{Dual of the modified principal sigma model.}

We set $\lambda=-1$ and $\mu=1$. Then
\[ E_e=K=-\h \sum_i f_i\tens f_i\]
Hence
\[ \hat {\bar E}^{-1}_{ij}=-\h(\delta_{ij}+\imath\eps_{ijk}\hat\pi_k),
\quad \hat\pi=2\vec t+\pmatrix{0\cr0\cr t^2}\]
and
\[ \hat {\bar E}_{ij}=-{2\over 1-t^2(t^2+4(t_3+1))}(\delta_{ij}
-\imath\eps_{ijk}\hat\pi_k-\hat\pi_i\hat\pi_j)\]
This defines the Lagrangian
\[ \hat{\cal L}={2\over 1-t^2(t^2+4(t_3+1))}\left(
\nabla_+\vec t\cdot\nabla_-\vec t-\imath\hat\pi\cdot(\nabla_+\vec t\times
\nabla_-\vec t)- (\hat\pi\cdot\nabla_+\vec t)(\hat\pi\cdot\nabla_-\vec t)
\right)\]
where $R_{t^{-1}*}t_\pm$ is computed as
\[ \nabla_\pm\vec t= \vec t_\pm-t_{\pm 3}{\vec t\over t_3+1}.\]
As before, the model in the form stated is complex
but with a different choice $\mu=\imath$ and different normalisation
of $r$ we can obtain a real model as well.

\subsubsection*{Dual of the pure-quasitriangular and principal sigma models.}

Here we set $\lambda=0$. Then rearranging $E_e$ above as an element of
$\gm\tens\gm$ we have
\[ E_e=-{1\over 4\mu}\left({1+2\mu\over 1+\mu}(f_1\tens f_1+f_2\tens f_2)
+{4\mu\over 1+2\mu}f_3\tens f_3+{\imath\over 1+\mu}f_1\wedge
f_2\right).\] One may then compute
\[ \hat {\bar E}_t=(E_e+\hat\Pi)^{-1}\]
and hence the Lagrangian. The result does not have any particular
simplifying features over the $\lambda=-1$ case above, so we omit
its detailed form.

Both limits of $\mu$ are singular, and require rescaling.
The $\mu\to\infty$ case makes sense after a rescaling of $r$
to $r/\mu$. This in turn scales the Lie cobracket of $\gg$ by $\mu^{-1}$
and hence also changes the Lie algebra structure of
$\gm$ to an Abelian one plus corrections of order $\mu^{-1}$. The effect of
this is to change the exponential map and the group law of $G^\star$, making
the latter Abelian. This can be expressed conveniently by working in
new coordinates with $\vec t$ scaled by $\mu^{-1}$. In this new
coordinate system we have
\[ \hat{\bar E}_t=-2\id+O(\mu^{-1})\]
since $\hat\Pi$ is linear in $\vec t$ to lowest order. The Lagrangian is
\[ \hat{\cal L}=2\vec t_+\cdot\vec t_-+O(\mu^{-1}).\]
Thus the dual model to the principal sigma model on $SU_2$ is an Abelian one
based on the group $\R^3$ with the usual linear wave equation.

The similar limit for the pure-quasitriangular case is ill-defined
since the Lie bracket of $\gm$ becomes singular as $\mu\to0$. Other scaling
limits of both the original model and its dual are possible in this case.

\note{A slightly different
limit is to scale $r$ as before by $\mu^{-1}$ and scale the Lie bracket
of $\gg$ by $\mu$ (so it becomes Abelian as $\mu\to 0$). This time the Lie
bracket of $\gm$ is unchanged but its Lie cobracket and hence $\hat\Pi$ becomes
scaled by $\mu$ and hence $O(\mu)$. }

\subsection{Point-particle limit of the quasitriangular models}

We have seen that the point-particle limit where $u$ is independent of
$x$ reduces to a classical mechanical dynamical system on the group $G$.
For our
quasitriangular models we have the following special cases.

\subsubsection*{Point-particle modified principal model.} From the
expressions for $E_u^{-1}$ etc. above, the Hamiltonian is
\eqn{conmodham}{ {\cal H}={1\over 4} K^{-1}((K\circ\Pi^R(u)+1)p,
(K\circ\Pi_2(u^{-1})+1)p)}
and the equations of motion are
\eqn{conmodeqm}{ u^{-1}\dot u=K^{-1}\circ((K\circ\Pi^R(u))^2-1)p,
\quad \dot p
=[K\circ\Pi^R(u)p,p].}
In this limit both the case entirely over $\R$ or the case where
$r$ and hence the cobracket are imaginary lead to well-defined real
equations of motion. In this case $\Pi$ is imaginary but so is the
Lie bracket of $\gm$ in the dual basis to the real basis of $\gg$.

For example, we can either work on $G=SL_2(\R)$ or, as more usual, on
$G=SU_2$. In the latter case (see above) we have
\[ K\circ\Pi^R(u)=\imath\pmatrix{0&-|b|^2&-\Im(a\bar b)\cr
|b|^2&0&\Re(a\bar b)\cr \Im(a\bar b)& -\Re(a\bar b)&0}=\imath\eps_{ijk}\pi_k
f_i\tens e_j.\]
Using the complexified Lie bracket on $su_2^\star$
we have the equations of motion for $p=p_if_i$ (with $p_i$ real) as
\[ \dot {\vec p}=\vec p (\vec p\times \pi)_3-p_3\vec p\times\pi\]
in terms of the vector cross product.
This can be written explicitly as
\[ \dot p_3=0,\quad \dot\rho=-{\imath\over 2}\bar a
b\rho^2+{\imath\over 2}a\bar b(|\rho|^2+2p_3^2)+\imath|b|^2\rho p_3,
\quad \rho\equiv p_1+\imath p_2\]
after a short computation. On the other hand,
\[ (K\circ\Pi^R)^2-1)_{ij}=(\pi^2-1)\delta_{ij}-\pi_i\pi_j\]
 hence the equation for $u$ is in our basis $e_i={-\imath\sigma_i\over 2}$ of
$su_2$ is
\[ u^{-1}\dot u=\imath(\pi^2-1)p\!\!\!/-\imath\pisl\pi\cdot p.\]
In our case $\pi^2=|b|^2$ and $\pi\cdot p=\Re(\rho \bar a b) -|b|^2p_3$, hence
\[ \dot u=-{\imath\over 2}\pmatrix{0&b\cr\bar b&0}(\rho\bar a b+\bar\rho a
\bar b-2|b|^2p_3)
-\imath |a|^2 u\pmatrix{p_3&\bar\rho\cr\rho&-p_3}.\]
Explicitly, this is
\[ \dot a=-\imath |a|^2(ap_3+b\rho),\quad \dot b=\imath b p_3
-{\imath\over 2}(1+|a|^2)a\bar \rho-{\imath\over 2}\rho\bar a b^2.\]
One may verify that this preserves $|a|^2+|b|^2=1$ as it must.

\subsubsection*{Point-particle pure-quasitriangular model.}
We set $\lambda=0$ and $E_u=E_e$ etc (the models are $G$-invariant). The
Hamiltonian and equations of motion are then
\eqn{conmuhamilt}{2 {\cal H}={1\over 1+2\mu}K((r_2+\mu K^{-1})p,
(r_2+\mu K^{-1})p)}
\eqn{conmueqm}{ u^{-1}\dot u=-{2\over 1+2\mu}(r_2+\mu K^{-1})\circ
K\circ(r_1+\mu K^{-1})p,
\quad \dot p={2\over 1+2\mu}[K\circ r_2p,p].}

Since these models are invariant, we know that $p\ra u^{-1}$ is conserved.
This means that we can let $Q=p(0)\ra u(0)^{-1}\in \gm$ be fixed and
substitute $p(t)=Q\ra u(t)$ into the equation for $\dot u$. We then solve
a first order non-linear differential equation for $u(t)$.

In particular, in the limit $\mu=0$ we obtain the $x$-independent
limit of the pure-quasitriangular model. Thus
\eqn{conpureeqm}{ {\cal H}={1\over 2}K(r_2p,r_2p),\quad u^{-1}\dot u
=2(r_2\circ K-1)r_2p,\quad
\dot p=2[K\circ r_2 p,p]}
using $r_1+r_2=K^{-1}$ to rearrange. In this case it makes sense to consider
the reduced variable $\xi=r_2p$ and write the equations of motion as
\eqn{xieqm}{ u^{-1}\dot u=2(r_2\circ K-1)\xi,\quad \dot\xi=2[r_2\circ
K \xi,\xi]} where we use that $r_2:\gg^\star\to\gg$ is a Lie
algebra homomorphism in view of the classical Yang-Baxter equation
(\ref{cybe})\cite{Ma:book}. We only need to solve this for $\xi$
in the image of $r_2$ but it is interesting that the equation
makes sense for any $\xi$ as an interesting integrable system on
the group manifold.

We can solve this for our strictly real form $\gg=sl_2(\R)$. We will solve
it here for the general (\ref{xieqm}); the special case of interest is
similar but more
elementary. Thus,
\[ r_2\circ K\pmatrix{H\cr X_+\cr X_-}=\pmatrix{\h H\cr X_+\cr 0}\]
so, writing $\xi(t)=h(t)H+x(t)X_++y(t)X_-$, we need to solve
\[ u^{-1}\dot u=-h(t)H-2y(t)X_-\]
and
\[ \dot h H+\dot x X_++\dot y X_-=[hH+2xX_+,hH+xX_++yX_-]
=2xy H-2hxX_+-2hyX_-,\]
which is the system of equations
\[ \dot h=2xy,\quad \dot x=-2hx,\quad \dot y=-2h y.\]
Note first of all that
\[{d\over d t}(h^2+\h \dot h)=0\]
so
\[ h^2+xy={\omega^2\over 4}\]
(say) is a constant. Inserting this into the equation for $h$ yields the
Riccati equation
\[ \dot h-\h \omega^2+2h^2=0\]
which has the general solution
\[ h(t)=\h \omega{\sinh(\omega t)+{2h(0)\over\omega}\cosh(\omega t)\over
\cosh(\omega t)+{2h(0)\over\omega}\sinh(\omega t)}.\]
We can then compute $y$ as
\[ y(t)=e^{-2\int^t_0 h(\tau)\extd \tau}y(0)\]
and similarly for $x(t)$. Since we only need $h,y$ to obtain $u(t)$ we can
consider the choice of $x(0)$ to be equivalent to the choice of $\omega$
(at least
in a certain range). The initial values of $h,y$ then determine their
general values
as above, and these then determine $u(t)$ given $u(0)$. The latter can
be expressed explicitly in terms of integrals on fixing a coordinate system for
$SL_2(\R)$.

For the point-particle limit of the pure quasitriangular models we are only
interested in $\xi\in b_+$ (the image of $r_2$), i.e. we specialise to
solutions
of the form $y(0)=0$, which clearly implies $y(t)=0$ and $\dot h=0$. In
this case
the solution is clearly
\[ \xi(t)={\omega\over 2}H+e^{-\omega t}x(0)X_+,\quad u(t)
=u(0)e^{-\h\omega t H}.\]
for initial data $\omega,x(0),u(0)$.

For the full physical momentum $p(t)$ we go back to (\ref{conpureeqm}).
If we write $p=2\omega\phi+x\psi_-+\bar x\psi_+$ say, then a similar
computation
using the Lie algebra of $sl_2(\R)^\star$ gives $\omega$ constant,
$\dot x=-\omega x$ as before, and additionally $\dot{\bar x}
=\omega \bar x$.
Hence the solution is
\[ p(t)=2\omega\phi+e^{-\omega t}x(0)\psi_-+e^{\omega t}\bar x(0)\psi_+,
\quad u(t)=u(0)e^{-\h\omega tH}\]
for constants $\omega,x(0),\bar x(0)$. As a check, it is easy to verify that
\[ Q_G=p\ra u^{-1}=(p\ra e^{\h\omega tH})\ra u(0)^{-1}\]
is conserved. Here $\psi_\pm\ra H=\mp2\psi_\pm$ and $\phi\ra H=0$ is the
relevant coadjoint action.

\subsubsection*{Point-particle principal model.} In the limit $\mu\to\infty$
of
(\ref{conmuhamilt})--(\ref{conmueqm}), we obtain the $x$-independent
limit of the principal sigma model. Here
\[ 4{\cal H}=K^{-1}(p,p),\quad u^{-1}\dot u=-K^{-1}\bar p,\quad
\dot{\bar p}=0\]
where $\bar p=\mu p$ is the renormalised momentum variable. This has the
general solution
\[ u(t)=u(0)e^{-tK^{-1}\bar p },\quad \bar p(t)=\bar p(0).\]
It is easy to see that $Q=\bar p\ra u^{-1}$ is constant as well,
using $K$ ad-invariant.

\section{Generalised T-Duality with double Neumann boundary conditions}

So far we have worked on providing a special class of Poisson-Lie T-dual
models within the established general framework. We now return to our
Hamiltonian
formulation of the general framework and observe that in this form the
main ideas can be extended to a much more general setting. Thus, from
 the symplectic form and the Hamiltonian we have just
calculated, we can see how the definition of T-duality could be
generalised. Begin with a Lie group $D$, with Lie algebra $\gd$,
and suppose that $\gd$
 is the direct sum of two
 subspaces $\cal E_-$ and $\cal E_+$. We take
 $\pi_+$ to be the projection to ${\cal E}_+$
with kernel ${\cal E}_-$, and $\pi_-$ to be the projection to ${\cal E}_-$
with kernel ${\cal E}_+$.

Suppose that there is a function $k:\Bbb
R^2\to D$, with the properties that $k_+k^{-1}(x_+,x_-)\in
\cal E_-$
and $k_-k^{-1}(x_+,x_-)\in \cal E_+$ for all $(x_+,x_-)\in \Bbb R^2$.
Then the relation $k_+k^{-1}(x_+,x_-)\in
\cal E_-$ can be summarised by $\pi_+(k_+k^{-1})=0$, and similarly
we get $\pi_-(k_-k^{-1})=0$. This gives the equations of motion
on
\[
\dot kk^{-1}\ =\ (\pi_--\pi_+)(k_xk^{-1})\ .
\]

Now we look at the symplectic form on the phase space.
Suppose that $\gd$ has an adjoint invariant inner product $\<,\>$.
If we imposed
boundary conditions that $k(0)$ and $k(\pi)$ were fixed, then
the symplectic form we computed earlier becomes
\[
2\omega(k;k_z,k_y)\ =\ \int^{\pi}_{x=0}
\<(k^{-1}k_y)_x,k^{-1}k_z\>\ dx\ .
\]
If we substitute $k_z=\dot t$, then we get
\[
2\omega(k;\dot k,k_y)\ =\ \int^{\pi}_{x=0}
\<k(k^{-1}k_y)_xk^{-1},\dot kk^{-1}\>\ dx\ =\
\int^{\pi}_{x=0}
\<(k_xk^{-1})_y,(\pi_--\pi_+)(k_xk^{-1})\>\ dx
\ .
\]
and so
\[
4\omega(k;\dot k,k_y)\ =\ -D_{(k;k_y)}
\int^{\pi}_{x=0}
\<k_xk^{-1},(\pi_+-\pi_-)(k_xk^{-1})\>\ dx
\ ,
\]
on the assumption that $\pi_+-\pi_-$ is Hermitian.
This will be true if the subspaces $\cal E_-$ and $\cal E_+$
are perpendicular with respect to the inner product.
Then we see that
$\omega(k;k_y,\dot k)=D_{(k;k_y)}{\cal H}(k)$, where
\[
4{\cal H}\ =\ \<  (\pi_{u+}-\pi_{u-})(u^{-1}u_x+s_xs^{-1}),
u^{-1}u_x+s_xs^{-1}\>\ .
\]
gives the Hamiltonian generating the time evolution.

The form of the boundary conditions we have imposed here should not come as
 too much of a surprise. Normally the string has boundary conditions
(for $k=us$ with $u\in G$ and $t\in M$)
$u_x=0$ at $x=0$ or $x=\pi$. This Neumann condition is designed to prevent
momentum
transfer out of the string at the edges. But if the system is to be
completely dual, we also need to impose a corresponding Neumann condition
on the
dual theory, which leads to the boundary condition $k_x=0$, the `double
Neumann' condition.
But then the equation of motion states $\dot k=0$ on the boundary.
Alternatively, if the reader prefers to work over $x\in\Bbb R$, we just
deal with
rapidly decreasing solutions. In either of these cases,
the symplectic form really is non-degenerate.

Now we have a phase space and Hamiltonian for the equations of
 motion
just based on an invariant inner product on $D$ and an orthogonal
 decomposition
$\cal E_-$ and $\cal E_+$ of $\gd$. If we take $D$ to be a
doublecross
product $D=G\dcross M$, and assume that the subspaces
${\rm Ad}_{u^{-1}}{\cal E}_\pm$
have graph coordinates $T_u$ and $E_u$ as before,
 we again recover the previous equations
of motion
for $u\in G$ in the factorisation $k=us$,
\[
\big(T_u(u^{-1}u_+)\big)_-\ -\ \big(E_u(u^{-1}u_-)\big)_+\ =\
\big[E_u(u^{-1}u_-),T_u(u^{-1}u_+)\big]\ .
\]
Importantly, we do not need to assume that the inner product has
any special properties
with respect to the decomposition $\gd=\gg+\gm$
 (such as being zero on $\gg$). We can also give the form of the Hamiltonian
for this general case:
\[
4{\cal H}\ =\ \<(E_u+I)(u^{-1}u_-),(E_u+I)(u^{-1}u_-)\>\ -\
\<(T_u+I)(u^{-1}u_+),(T_u+I)(u^{-1}u_+)\>\ .
\]
The corresponding dual formula would produce exactly the same value.

\subsection{Poisson brackets and the central extension}
In this section we continue with the generalised T-duality
and boundary conditions of the last section.
The phase space for our system is infinite dimensional, so it is rather hard to
describe the functions on it directly. We shall describe a `nice' set of
functions,
and hope that more general functions are expressible as a product of these
nice
functions.

If $v\in C^\infty((0,\pi),\gd)$, we can look at the vector field
$k_z=vk$ for $k\in  C^\infty((0,\pi),D)$. To preserve the boundary
conditions we consider only those
$v\in C^\infty((0,\pi),\gd)$ which tend to zero at the end points.
Consider
\[
\omega(k;k_y,k_z)\ =\ -\frac12 \int\<(k^{-1}k_y)_x,k^{-1}k_z\>\, dx\ =\
-\frac12 \int\<(k_x k^{-1})_y,v\>\, dx\ =\
-\frac12 D_{(k;y)} \int\<k_x k^{-1},v\>\, dx\ .
\]
It follows that the function which acts as a Hamiltonian generating this flow
is
\[
f_v(k)\ =\ -\frac12  \int\<k_x k^{-1},v\>\, dx\ .
\]
We can calculate the Poisson brackets between these nice functions
quite easily:
\[
\{f_v,f_w\}\ =\ f_v'(k,wk)\ =\ f_{[v,w]}\ -\ \frac12\int \< w_x,v\>\, dx\ .
\]
We now see the appearance of a central extension term in the Lie algebra.
The Poisson brackets can be written as $
\{f_v,f_w\}= f_{[v,w]}+\vartheta(v,w)f_c$,
where $f_c(k)=1$ and the cocycle $\vartheta(v,w)=-\int \< w_x,v\>\, dx/2$.
We can also manufacture a derivation term, which corresponds to the
momentum (the operation of incrementing the $x$ coordinate). Consider
\begin{eqnarray*}
\omega(k;k_y,k_x)&=& -\frac12 \int\<(k^{-1}k_y)_x,k^{-1}k_x\>\, dx\ =\
-\frac12 \int\<(k_x k^{-1})_y,k_x k^{-1}\>\, dx \cr &=&
-\frac14 D_{(k;y)} \int\<k_x k^{-1},k_x k^{-1}\>\, dx\ .
\end{eqnarray*}
Thus the momentum is given by
\[
f_d(k)\ =\ -\frac14  \int\<k_x k^{-1},k_x k^{-1}\>\, dx\ .
\]
A brief calculation shows that $\{f_d,f_v\}=f_{v'}$ and $\{f_d,f_c\}=0$.

\subsection{Adjoint symmetries of the model and dual model}
In this section we consider the left multiplication symmetry again, however
this time we
can simultaneously describe the action on the dual models. This requires
some care
with the boundary conditions, and we shall take the double Neumann
condition on loops,
i.e.\ $k=e$ and $k_x=0$ at both boundaries. The operation of left
multiplication
by constants does not preserve these conditions, but we can use our freedom
to introduce a right multiplication to work with the adjoint action instead.

Take the action on the phase space given by ${\rm Ad}_d$ for $d\in D$. This
preserves the
boundary conditions, and preserves the models in the case where ${\rm
Ad}_d{\cal E}_\pm=
{\cal E}_\pm$. The corresponding infinitesimal motions are generated by the
moment map
\[ I_\delta(k)=-\h\int \<k_x k^{-1},\delta\> \extd x\ ,\quad \delta\in\gd.\]
If the map ${\rm ad}_\delta$ preserves the subspaces ${\cal E}_\pm$ then
this formula
gives conserved charges for the system.

\subsection{Automorphism symmetries of the model and dual model}
Here we consider symmetries of the phase space arising from
group automorphisms $\theta:D=G\bowtie M\to D$. This is really a generalisation
of the previous subsection, where we just considered automorphisms given by
the adjoint action, i.e.\ inner automorphisms. We consider the same
boundary conditions as in the last subsection. For convenience we also
assume that
the two subspaces $\theta{\cal E}_\pm$ of the Lie algebra $\gd$
are perpendicular for the given inner product. This is not really needed,
as we can
always manufacture
a new ${\rm Ad}$-invariant inner product from the old one using the
automorphism
in order to make this true.

Given these conditions, any automorphism $\theta:D\to D$ will induce a map
$\tilde\theta$ on the phase space given by $(\tilde\theta
k)(x)=\theta(k(x))$. This map will be
symplectic if $\theta$ preserves the given inner product on $\gd$, and if
$\theta {\cal E}_\pm={\cal E}_\pm$ then the map will preserve the given
models.
In general $\tilde\theta k$ will factor to give $G$-models and dual $M$-models
which are a mixture of the original $G$-models and dual $M$-models
given by factoring $k$. However there are two special cases
worthy of mention.

1)\quad
The automorphism $\theta:D\to D$ is called subgroup preserving if $\theta
G\subset G$
and $\theta M\subset M$. In this case a factorisation $k=us$ for $u\in G$
and $s\in M$
is sent to $\theta(k)=\theta(u)\theta(s)$, and $\theta(u)$ is a solution
of the sigma model on $G$. In the same manner, if $t$ is a solution of the
sigma model on
$M$, then $\theta(t)$ is also a solution
of the sigma model on $M$.

2)\quad
The automorphism $\theta:D\to D$ is called subgroup reversing if $\theta
G\subset M$
and $\theta M\subset G$. If such an automorphism exists, the double
$D=G\bowtie M$
is called self-dual\cite{BegMa:dif}.  In this case a factorisation $k=us$ 
for $u\in G$ and $s\in M$
is sent to $\theta(k)=\theta(u)\theta(s)$, and $\theta(u)$ is a solution
of the sigma model on $M$. In the same manner, if $t$ is a solution of the
dual sigma model on
$M$, then $\theta(t)$ is also a solution
of the sigma model on $G$. In this manner the solutions of the sigma model
on $G$
and the dual sigma model on $M$ are related by a group homomorphism
from $G$ to $M$, and in that sense the models are self-dual.

Other symmetries may be constructed. For example of we have $\theta{\cal
E}_+={\cal E}_-$
and $\theta{\cal E}_-={\cal E}_+$ then the map $\hat
\theta(k)(t,x)=\theta(k(t,\pi-x))$
sends a solution $k$ of the model into another solution.
\bigskip

The explicit computation of examples of our generalised T-duality
along the above lines is a topic for further work. However, the data
required for the construction do exist in abundance. For example, given any
two Lie algebras $\gg_0\subset\gd$ whose Dynkin diagrams differ by the
deletion of some
nodes, one has an inductive construction $\gd=(\gn\lcross\gg_0)\dcross\gn^*$
where $\gn$ are braided-Lie bialgebras \cite{Ma:blie}. For a concrete example,
one has, locally,
\[ D=SO(1,n+1)=(\R^n\lcross SO(n))\dcross\R^n\]
as the decomposition of conformal transformations into Poincar\'e and
special conformal translations. The group $D$ has a non-degenerate bilinear 
form as required (although not positive-definite). The explicit construction 
of the required factorisation and the associated bicrossproduct quantum groups
and T-dual models will be attempted elsewhere.

%\bibliographystyle{unsrt}
%\bibliography{biblio}

\end{document}